\documentclass[11pt]{article}
\usepackage{latexsym}
\usepackage{cite}
\usepackage{amsmath}           
\usepackage{amsfonts}          
\usepackage{amssymb}           
\usepackage{bm}                
\usepackage[english]{babel}

\usepackage{hyperref}

\newcommand{\be}{\begin{equation}}
\newcommand{\ee}{\end{equation}}
\newcommand{\ef}[1]{\, #1} 

\newcommand{\bydef}{:=}
\def\reff#1{(\protect\ref{#1})}
\newcommand{\sgn}{\mathop{\rm sgn}\nolimits}
\newcommand{\diag}{\mathop{\rm diag}\nolimits}
\newcommand{\rowdet}{\mathop{\hbox{\rm row-det}}\nolimits}
\newcommand{\coldet}{\mathop{\hbox{\rm col-det}}\nolimits}
\newcommand{\symdet}{\mathop{\hbox{\rm sym-det}}\nolimits}

\newcommand{\scrs}{{\mathcal{S}}}
\newcommand{\scra}{{\mathcal{A}}}

\def\R{{\mathbb R}}
\def\C{{\mathbb C}}
\def\N{{\mathbb N}}
\def\Q{{\mathbb Q}}

\newcommand{\psibar}{\bar{\psi}}
\newcommand{\T}{{\rm T}}

\newcommand{\aaa}{a}
\newcommand{\adg}{a^{\dagger}}
\newcommand{\ket}[1]{\left| #1 \right\rangle}
\newcommand{\bra}[1]{\left\langle #1 \right|}
\newcommand{\braket}[2]{\left\langle #1 | #2 \right\rangle}

\newcommand{\nunu}[1]{( #1 )}

\newtheorem{defin}{Definition}[section]
\newtheorem{prop}[defin]{Proposition}
\newtheorem{lemma}[defin]{Lemma}

\newtheorem{coroll}[defin]{Corollary}

\def\proof{\par\medskip\noindent{\sc Proof.\ }}
\newcommand{\qed}{\hfill $\Box$\vspace{1.5mm}}

\def\proofof#1{\bigskip\noindent{\sc Proof of #1.\ }}

\begin{document}

\title{ Noncommutative determinants, \\
           Cauchy--Binet formulae, \\
           and Capelli-type identities  \\[5mm]
 \large\bf II.~Grassmann and quantum oscillator algebra representation \\[1cm]}

\author{
  { Sergio Caracciolo}                       \\
  {\small\it Dipartimento di Fisica dell'Universit\`a degli Studi di Milano and INFN}      \\
  {\small\it via Celoria 16,}                       
  {\small\it I-20133 Milano, ITALY}                \\
  {\small\tt Sergio.Caracciolo@mi.infn.it}         \\
  {\protect\makebox[5in]{\quad}}
    \\
  { Andrea Sportiello\thanks{On leave from
         Dipartimento di Fisica and INFN, Universit\`a degli Studi di Milano,
         via Celoria 16, I-20133 Milano, ITALY.}}  \\
  {\small\it LIPN, UMR CNRS 7030}                  
  {\small\it Universit\'e Paris-Nord}              \\
  {\small\it 99, avenue Jean-Baptiste Cl\'ement}   
  {\small\it 93430 Villetaneuse, FRANCE}           \\
  {\small\tt Andrea.Sportiello@lipn.fr}            \\
  {\protect\makebox[5in]{\quad}}
\\[1cm]
\\
{\small
05A19, 05A30, 05E15, 13A50,
15A54,
16T20,
20G42
(Secondary).}
}
\maketitle
\thispagestyle{empty}   




\newpage
\vspace*{4cm}

\begin{abstract}
We prove that, for $X$, $Y$, $A$ and $B$ matrices with entries in a
non-commutative ring such that
\hbox{$[X_{ij},Y_{k\ell}]=-A_{i\ell} B_{kj}$},
satisfying suitable commutation relations
(in particular, $X$ is a Manin matrix),
the following identity holds
\begin{equation*}
\coldet X \coldet Y 
= 
\bra{0}
\coldet (\aaa A + X (I-\adg B)^{-1} Y)
\ket{0}
\ef.
\end{equation*}
Furthermore, if also $Y$ is a Manin matrix,
\begin{equation*}
\coldet X 
\coldet Y
=\int \mathcal{D}(\psi, \psibar)
\exp 
\Bigg(
\sum_{k \geq 0}
\frac{(\psibar A \psi)^{k}}{k+1}
(\psibar X B^k Y \psi)
\Bigg)
\ef.
\end{equation*}
Notations: $\bra{0}$, $\ket{0}$, are respectively the bra and the ket
of the ground state, $\adg$ and $\aaa$ the creation and annihilation
operators of a quantum harmonic oscillator, while $\psibar_i$ and
$\psi_i$ are Grassmann variables in a Berezin integral.  These results
should be seen as a generalization of the classical Cauchy--Binet
formula, in which $A$ and $B$ are null matrices, and of the
non-commutative generalization, the Capelli identity, in which $A$ and
$B$ are identity matrices and
$[X_{ij},X_{k\ell}]=[Y_{ij},Y_{k\ell}]=0$.
\end{abstract}



\newpage


\section{Introduction}

\subsection{The Cauchy--Binet theorem}

Let $R$ be a commutative ring, and let $M = (M_{ij})_{i,j=1}^n$
be a $n \times n$ matrix with elements in $R$.
The determinant of the matrix $M$ can be defined as
\begin{equation}
   \det M  \;\bydef\;
   \sum_{\sigma \in \scrs_n} \sgn(\sigma) \, 
   M_{\sigma(1) 1} \, M_{\sigma(2) 2} \,\cdots\, M_{\sigma(n) n}
   \;
 \label{def.det}
\end{equation}
where $\scrs_n$ is the permutation group of the set
$[n]=\{1,2,\dots,n\}$  and $\sgn(\sigma)$ is the sign of the
permutation $\sigma$.

Let $X$ be a $n \times m$ matrix and $Y$ a $m \times n$ matrix with
elements in the commutative ring $R$.  For each subset $I\subseteq
[m]$ let be $X_{[n],I}$ the minor of $X$ with columns in $I$ and
similarly $Y_{I,[n]}$ the minor of $Y$ with rows in $I$. The classical
Cauchy--Binet formula relates the product of the determinant of these
matrices to the determinant of the product. More precisely 
\be
\label{form.CB}
\sum_{\substack{L \subseteq [m] \\  |L|= n}} \det X_{[n],L} \det
Y_{L,[n]} 
=  \det (XY) \, . 
\ee
In order to generalize the definition \reff{def.det} to matrices with
elements in a {\em noncommutative}\/ ring $R$, the first problem
encountered is that it is ambiguous without an ordering prescription
for the product.  Rather, numerous alternative ``determinants'' can be
defined: for instance, the {\em column-determinant}\/
\begin{equation}
   \coldet M  \;\bydef\;
   \sum_{\sigma \in \scrs_n} \sgn(\sigma) \, \prod_{i=1}^n M_{\sigma(i)\, i} 
 \label{def.coldet}
\end{equation}
and the {\em row-determinant}\/
\begin{equation}
   \rowdet M \;\bydef\; 
   \sum_{\sigma \in \scrs_n} \sgn(\sigma) \,
   \prod_{i=1}^n M_{i \, \sigma(i) }
 \label{def.rowdet}
 \;.
\end{equation}
(Note that $\coldet M = \rowdet M^{\T}$.)  It is intended above that,
when dealing with non-commuting quantities having indices depending on
a single integer, the product symbol $\Pi$ denotes
an ``ordered product'', i.e.
\be
\prod_{i=k}^{k+\ell} f_i := f_k f_{k+1} \dots f_{k+\ell}
\ef.
\ee
In~\cite{uscapelli} we have proven, in collaboration with
A.\,D.~Sokal, non-commutative generalizations of the Cauchy--Binet
formula.
In order to express our result, we called the matrix $M$
{\em column-pseudo-commutative} in the case
\be
   [M_{ij}, M_{k\ell}]  = [M_{i\ell}, M_{kj}]  \qquad\hbox{\rm for all } i,j,k,\ell
 \label{def.colpc.1}
\ee
and
\be
   [M_{ij}, M_{i\ell}]  =  0   \qquad\hbox{\rm for all } i,j,\ell
   \;.
 \label{def.colpc.2}
\ee
(Similarly, we said a matrix $M$ to be {\em row-pseudo-commutative}\/ in
case $M^\T$ is column-pseudo-commutative)\footnote{Note that
  (\ref{def.colpc.1}) implies $2 [M_{ij}, M_{i\ell}] = 0$,
  i.e.\ \emph{twice} equation (\ref{def.colpc.2}), a subtlety, of
  relevance only when the field $K$ over which the ring $R$ is defined
  is of characteristic $2$, that will appear several times along the
  paper.}.  Furthermore, we said that $M$ has
\emph{weakly column-symmetric (and row-antisymmetric) commutators}
if (\ref{def.colpc.1}) holds for $i \neq k$
(and (\ref{def.colpc.2}) not necessarily holds).


We proved~\cite[Proposition 1.2]{uscapelli} that\footnote{Here
  we perform a change of notation for future convenience $(A^\T \to X,
  B \to Y, h \to A)$ and consider only the case $r=n$.}
\begin{prop}[noncommutative Cauchy--Binet] 
\label{theo.capCB}
Let $R$ be a
ring, and let $X$ be a $n \times m$ matrix and $Y$ a $m \times n$
matrix with elements in $R$.  Suppose that
\be
  [X_{ij},Y_{k\ell }] =  \, -A_{i\ell } \delta_{kj}  \qquad\hbox{\rm for all } i,j,k,\ell
\ee
with $A$ a $n \times n$ matrix.
Then 
\begin{itemize}
\item[(a)] If $X$ is row-pseudo-commutative, then
\be
\label{eq.thCapCB}
\sum_{\substack{L \subseteq [m] \\  |L|= n}} \coldet X_{[n],L} \coldet
Y_{L,[n]} = \coldet ( XY +Q^{\rm col} )
\ee
where
\be
\label{eq.defQcol}
Q^{\rm col}_{ij} \bydef A_{i j} (n-j) \, .
\ee
\item[(b)] If $Y$ is column-pseudo-commutative, then
\be
\sum_{\substack{L \subseteq [m] \\  |L|= n}} \rowdet X_{[n],L} \rowdet
Y_{L,[n]} = \rowdet ( XY + Q^{\rm row} )
\ee
where
\be
\label{eq.defQrow}
Q^{\rm row}_{ij} \bydef A_{i j} (i-1) \, .
\ee
\item[(c)] In particular, if $[X_{ji},X_{\ell k}] = 0$ and
  $[Y_{ij},Y_{k\ell}] =0$ whenever $j\neq \ell$, then
\be
\sum_{\substack{L \subseteq [m] \\  |L|= n}} \det X_{[n],L} \det
Y_{L,[n]} = \coldet ( XY +Q^{\rm col} ) = \rowdet ( XY + Q^{\rm row}
)\, .
\ee
\end{itemize}
\end{prop}
\noindent
With respect to the commutative case~\reff{form.CB}, the determinants
are replaced by one of its non-commutative generalizations, but the
left-hand side keeps the same form, while on the right-hand side the
product $XY$ requires an additive correction.

An example of a non-commutative ring $R$ is the Weyl algebra $A_{m \times n}(K)$
over some field $K$ of characteristic 0 (e.g.\ $\Q$, $\R$ or $\C$)
generated by a $m \times n$ collection $Z = (z_{ij})$
of commuting indeterminates (``positions'')
and the corresponding collection $\partial = (\partial/\partial z_{ij})$
of differential operators (proportional to ``momenta''); so that
\begin{subequations}
\begin{align}
\left[z_{ij}\, , \frac{\partial}{\partial z_{k \ell }}\right]
= & \,
-\delta_{i k} \delta_{j \ell}
\ef;
\\ 
[z_{ij}\, , z_{k\ell }] =  & \, \left[\frac{\partial}{\partial
    z_{ij}}\, ,  \frac{\partial}{\partial z_{k \ell }}\right] = 0
\ef.
\end{align}
\label{eqs.weyldefs}%
\end{subequations}
If we set $m=n$, $X=Z^\T$ and $Y=\partial$, we soon get $A_{ij} =
\delta_{ij}$ for each $i,j \in [n]$ and
\begin{align}
\det X \det \partial = & \coldet [ X^\T \partial + \diag (n-1, n-2,
  \dots, 0) ]  \label{Capelli:col}\\
= & \rowdet [  X^\T \partial + \diag (0, 1, \dots, n-1) ]
\label{Capelli:row}
\end{align}
which are the Capelli identities
\cite{Capelli_1882, Capelli_1887, Capelli_1888, Capelli_1890} of
classical invariant theory \cite{Weyl_46, Procesi, Howe_89}, a field
of research that, in more than a century, has remained active up to
recent days (a forcerly incomplete selection of papers on the subject
includes
\cite{Turnbull_48, Wallace_53, Kostant_91, Kostant_93, Howe_91,
  Itoh_00, Itoh-Umeda_01, Noumi_94, Noumi_96, Umeda_98, Kinoshita_02,
  Mukhin_06, Molev_99, Nazarov_97, Nazarov_98, Okounkov_96b}).
Because of this example, the correction term due to the presence of
the matrix $Q$ which appears in the non-commutative case is sometimes
called the ``quantum'' correction with respect to the formula in the
commutative case~\reff{form.CB}.

Chervov, Falqui and Rubtsov give in~\cite{Chervov_09} an extremely
interesting survey of the algebraic properties of
row-pseudo-commutative matrices (which they call ``Manin matrices'',
because a similar notion has proven fruitful in the context of quantum
groups, where it arose already two decades ago in Manin's work
\cite{Manin_87,Manin_88,Manin_91,Chervov_08}), when the ring
$R$ is an associative algebra over a field of characteristic $\neq
2$. In particular, \cite[Section~6]{Chervov_09} contains an
interesting generalization of our result.  Another recent interesting
survey, on combinatorial methods in the study of non-commutative
determinants, is the PhD Thesis of M.~Konvalinka~\cite{KonvPhd}.

In this paper we will investigate a stronger version of
Proposition~\ref{theo.capCB}. In particular we relax the condition
that for all $i,j,k,\ell$
\be
  [X_{ij},Y_{k\ell }] =  -A_{i\ell } \delta_{kj}  \label{old}
\ee
to
\be
  [X_{ij},Y_{k\ell }] =  -A_{i\ell } B_{kj}  \label{new}
\ee
where $B$ is a $m \times m$ matrix whose elements are supposed to
commute with everything.

Remark that, whenever $B$ is invertible,\footnote{Recall that, in our
  case, this is not just a matter of the matrix being non-singular: as
  the entries $B_{ij}$ are valued in a ring, not even the
  single entries, even when non-zero, are guaranteed to have a
  multiplicative inverse, i.e.\ not even the case $n=1$ is easy.}
from~\reff{new} by
multiplication of $B_{js}^{-1}$ and sum over $j$ we get
\be
[ (XB^{-1})_{is}, Y_{k\ell }] = -A_{i\ell } \delta_{ks}
\ee
which is of the form~\reff{old}, and similarly by multiplication of
$B_{sk}^{-1}$ and sum over $k$
\be
[ X_{ij}, (B^{-1}Y)_{s\ell }] = -A_{i\ell } \delta_{sj}
\label{B-1commrel}
\ee
and, as if $X$ is row-pseudo-commutative also $XB^{-1}$ is such, while
if $Y$ is column-pseudo-commutative also $B^{-1}Y$ is such. Thus,
quite trivially, Proposition~\ref{theo.capCB} can be used to express,
for example in the case (a)
\be
\sum_{\substack{I, L \subseteq [m] \\ |L|= |I| = n}} \coldet
X_{[n],I}\, \det B^{-1}_{IL}\, \coldet Y_{L,[n]} = \coldet ( XB^{-1}Y
+Q^{\rm col} )\, . 
\label{soft}
\ee
In agreement with the philosophy of the original Capelli identity, our
goal in this paper is in another direction: we want to find
generalizations of Proposition \ref{theo.capCB}, under the more
general (\ref{new}), in which the left-hand side of (\ref{eq.thCapCB})
(and variants) is kept \emph{exactly in this form} (with no dependence
from $B$ whatsoever), and investigate for a generalized ``quantum
correction'' on the right-hand side.

We have not been able to reach an expression as simple as we got
previously in Proposition \ref{theo.capCB},
(not even in the case when $B$ is invertible).
However, we have found closed formulas with the help of the algebra
and the Hilbert space of a single ``bosonic quantum oscillator'' (also
known as \emph{Heisenberg-Weyl Algebra}),
and, also, as a Berezin integral in Grassmann algebra, corresponding to
``fermionic quantum oscillators'' (see respectively the following
Propositions \ref{new.capCB} and~\ref{prop:grass}, which are the main
results of the paper).

We point out here a possible source of confusion. While, at the
foundations of invariant theory, Capelli identities have been
discovered within their explicit realization in Weyl Algebra (the
example of equations (\ref{eqs.weyldefs})), it is nowadays clear, and
along the lines e.g.\ of \cite{uscapelli}, \cite{Chervov_09}, and
several other papers, that the appropriate context of this family of
identities is the identification of sufficient conditions on the
commutation rules for the elements of the involved matrices,
regardless from the presentation of rings $R$, and matrices valued in
$R$, realizing these rules. To characterize and classify these
realizations (or, even, to determine their existence) is a problem
that we find important, but of separate interest, and we do not treat
it here. The role of the Weyl-Heisenberg and Grassmann algebras
mentioned above is \emph{not} at the level of the explicit realization
of the matrices. It consists instead of an auxiliary structure,
implementing certain combinatorial relations at the level of
manipulation of commutators, that arise along the lines of the proof.

We annotate here an interesting paper, by Blasiak and
Flajolet~\cite{Flajolet}, presenting a collection of classical and new
facts on the role of
Weyl-Heinsenberg Algebra in combinatorics, in the spirit of the
discussion above.

\subsection{The bosonic quantum oscillator}
\label{ssec.bosointro}

Following the classical treatment of the quantum oscillator by
Dirac~\cite[Chapter 6]{Dirac}, let us introduce the operator $a$ and
its adjoint $\adg$, called respectively {\em annihilation} and
{\em creation} operator, and the Hermitian number operator $N =
\adg \aaa$.

They satisfy the commutation relations of the Weyl-Heisenberg algebra
\begin{align}
[\aaa, \adg] 
&= 1
\ef,
&
[N, \aaa] 
&= -\aaa
\ef,
&
[N, \adg] 
&= \adg
\ef.
\label{WH}
\end{align}
Let $\ket{n}$ with $n\in \N$ be the eigenstate of $N$ corresponding
to the eigenvalue $n$, that is
\be
N \ket{n}
= n \ket{n}
\ef.
\ee
In particular the lowest eigenstate of $N$, $\ket{0}$, is annihilated
by $\aaa$
\be
\aaa
\ket{0} 
= 0
\ef.
\ee
Without loss of generality, we assume it to be of unit norm,
$\braket{0}{0}=1$.

Our first generalization of the Capelli identity is stated within this framework.
\begin{prop}
\label{new.capCB}
Let $R$ be a 
ring,
and let $X$ be a $n \times m$ matrix and $Y$ a $m \times n$ matrix with elements in $R$. 
Suppose that
\be
\label{XY}
  [X_{ij},Y_{k\ell }] =  \, -A_{i\ell } B_{kj}  \qquad\hbox{\rm for all } i,j,k,\ell
\ee
with $A$ a $n \times n$, and $B$ a $m \times m$ matrix whose elements commute with everything.
Then 
\begin{itemize}
\item[(a)] If $X$ is row-pseudo-commutative, 
and
\be
[X_{ij}, A_{k\ell}] - [X_{kj}, A_{i\ell}] = 0 
\qquad\hbox{\rm for all } i,j,k,\ell 
\label{XA2}
\ee
then
\be
\sum_{\substack{L \subseteq [m] \\ |L|= n}} \coldet X_{[n],L} \coldet
Y_{L,[n]} = \langle 0 | \coldet [ a\,  A + X (1-\adg\, B)^{-1}Y ]
\, | 0 \rangle\, .
\ee
\item[(b)] If $Y$ is column-pseudo-commutative, 
and
\be
[Y_{ij}, A_{k\ell}] - [Y_{i\ell}, A_{kj}] = 0
\qquad\hbox{\rm for all } i,j,k,\ell
\label{XA2caseb}
\ee
then
\be
\sum_{\substack{L \subseteq [m] \\ |L|= n}} \rowdet X_{[n],L} \rowdet
Y_{L,[n]} = \langle 0 | \rowdet [ \adg\,  A + X (1-a\, B)^{-1}Y ]
\, | 0 \rangle\, .
\ee
\item[(c)] In particular, if $[X_{ji},X_{\ell k}] = 0$ and
  $[Y_{ij},Y_{k\ell}] =0$ whenever $j\neq \ell$,
then
\begin{align}
\sum_{\substack{L \subseteq [m] \\ |L|= n}} \det X_{[n],L} \det
Y_{L,[n]}= & \langle 0 | \coldet [ a\,  A + X (1-\adg\,
  B)^{-1}Y ] \, | 0 \rangle \\
= & \langle 0 | \rowdet [ \adg\,  A + X (1-a\, B)^{-1}Y ]
\, | 0 \rangle\,
\, .
\end{align}
\end{itemize}
\end{prop}
The further commutation condition (\ref{XA2}) (and the counterpart
(\ref{XA2caseb}) for case (b)) appears as a subtle technicality, that
we did not succeed to avoid. Note however that, as shown in Lemmas
\ref{lem.XABcond} and \ref{lem.XABcond2} through an analysis of the
consequences of the Jacobi Identity, it is implied by a very mild
condition on $B$, (informally, that two vectors $\vec{u} , \vec{v} \in
R^m$ exist such that the scalar product $(\vec{u}, B \vec{v})$ is a
regular element of the ring, i.e., it is not zero, and not a divisor
of zero). In particular, this is obviously the case under the
circumstances originally treated in \cite{uscapelli}, where $B=I$.

%

As an example, let the non-commutative ring $R$ be the Weyl algebra $A_{m \times s}(K)$
over some field $K$ of characteristic 0 (e.g.\ $\Q$, $\R$ or $\C$)
generated by an $m \times s$ collection $Z = (z_{i}^a)$ with $i\in [n]$ and $a\in [s]$
of commuting indeterminates 
and the corresponding collection $\partial = (\partial/\partial z_{i}^a)$
of differential operators; so that
\begin{subequations}
\label{eq.commuteasy}
\begin{align}
\bigg[z_{i}^a\, , \frac{\partial}{\partial z_{j}^b}\bigg]
&=
-\delta_{i j} \delta^{a b}
\ef;
\\ 
[z_{i}^a\, , z_{j}^b] =
\bigg[
\frac{\partial}{\partial z_{i}^a}\, ,
\frac{\partial}{\partial z_{j}^b}\bigg] 
&= 0
\ef.
\end{align}
\end{subequations}
Let 
\begin{align}
X_{ij}
&=
\sum_{a=1}^s z_i^a \alpha_j^a
\ef,
&
Y_{k \ell}
&=
\sum_{a=1}^s \beta_k^a \frac{\partial}{\partial z_\ell^a}
\ef,
\end{align}
with $\alpha_j^a, \beta_k^a$ commuting with everything, so that for
all $i,\ell\in [n]$ and $j , k \in [m]$
\be
[X_{ij}, X_{k \ell}] = [Y_{ij}, Y_{k \ell}] = 0 
\ee
and
\be
[X_{ij}, Y_{k \ell}] = - \delta_{i \ell}  \sum_{a=1}^s \beta_k^a \alpha_j^a
\ee
which, in our notation means that
\begin{align}
A_{i \ell} 
&= \delta_{i \ell}
\ef,
&
B_{k j} 
&= \sum_{a=1}^s \beta_k^a \alpha_j^a
\ef.
\end{align}
Remark that the rank of the $m \times m$ matrix $B$ is $\min(m,s)$,
in particular, when $s < m$, $B$ is not invertible.

In the particular case in which $B_{ij} = \delta_{ij}$ for each
$i,j\in [m]$, both Proposition~\ref{theo.capCB} and~\ref{new.capCB}
apply. As a consequence, the right hand sides must be equal and, for
example, if $X$ is row-pseudo-commutative, then
\be
\coldet ( XY +Q^{\rm col} ) 
= \langle 0 | \coldet [ a\,  A +
  (1-\adg)^{-1} XY ] \, | 0 \rangle
\ee
while, if $Y$ is column-pseudo-commutative, then
\be
\rowdet ( XY +Q^{\rm row} ) 
= \langle 0 | \rowdet [ \adg\,  A +
  (1-a)^{-1} XY ] \, | 0 \rangle \, .
\ee
These relations are indeed valid regardless from the fact that $A$ is
related to the commutator of $X$ and $Y$, i.e.\ they are a consequence
of a stronger fact
\begin{prop}
\label{prop.old}
Let $R$ be a ring and $U$ and $V$ be two $n\times n$ matrices with
elements in $R$. Then
\be
\coldet (U + Q^{\rm col}) 
=
\bra{0}
\coldet \big( a V + (1-\adg)^{-1} U \big)
\ket{0}
\ee
where
\be
Q^{\rm col}_{ij} \bydef V_{i j} (n-j)
\ef,
\ee
and
\be
\rowdet (U + Q^{\rm row})
=
\bra{0}
\rowdet \big( \adg V + (1-a)^{-1} U \big)
\ket{0}
\label{old.row}
\ee
where
\be
Q^{\rm row}_{ij} \bydef V_{i j} (i-1)
\ef.
\ee
\end{prop}
This fact, together with a generalization, is proven in Section \ref{sec.equiv}.


\subsection{The Grassmann algebra}

The determinant of a $n\times n$ matrix $M$ with elements in a
commutative ring can be represented as a Berezin integral over the
Grassman algebra generated by the $2n$ anti-commuting variables
$\{\psi_i, \psibar_i\}_{i\in [n]}$ (for an introduction to such a
topic we invite the interested reader to refer to~\cite[Appendix
  B]{uscayley}). More precisely: 
\be
\det M = \int \mathcal{D}(\psi, \psibar)
\exp ( \psibar M \psi ) \label{simplegrass}
\ef
\ee
where
\be
\mathcal{D}(\psi, \psibar) :=  \prod_{i=1}^n d\psi_i\,d\psibar_i
\ef.
\ee
Therefore the Cauchy--Binet theorem can also be written as the identity
\be
\sum_{\substack{L \subseteq [m] \\ |L|= n}} \det X_{[n],L} \det
Y_{L,[n]} = \int \mathcal{D}(\psi, \psibar) \exp ( \psibar XY \psi
)\, . \label{CBclassico}
\ee
We have obtained the following generalization
\begin{prop}
\label{prop:grass}
Let $R$ be a 
ring containing the rationals, and let $X$ be a $n \times m$ matrix
and $Y$ a $m \times n$
matrix with elements in $R$. Suppose that
\be
  [X_{ij},Y_{k\ell }] =  \, -A_{i\ell } B_{kj}  \qquad\hbox{\rm for all } i,j,k,\ell \label{hyp1}
\ee
with $A$ a $n \times n$, and $B$ a $m \times m$ matrix whose elements commute with everything.
Let $I_m$ the $m\times m$ identity matrix. 
Assume that
\begin{align}
\label{eq.XAgrass}
[X_{ij}, A_{k\ell}] - [X_{kj}, A_{i\ell}] 
&= 0
\qquad\hbox{\rm for all } i,j,k,\ell
\ef;
\\
\label{eq.YAgrass}
[Y_{ij}, A_{k\ell}] - [Y_{i\ell}, A_{kj}] 
&= 0
\qquad\hbox{\rm for all } i,j,k,\ell
\ef.
\end{align}
Then 
\begin{itemize}
\item[(a)] If $X$ and $Y$ are row-pseudo-commutative,
then
\be
\label{eq.inPropGrass}
\begin{split}
&
\hspace{-10mm}
\sum_{\substack{L \subseteq [m] \\ |L|= n}} 
\coldet X_{[n],L} \coldet Y_{L,[n]}
=
\int \mathcal{D}(\psi, \psibar)
\exp \bigg(
\sum_{k \geq 0}
\frac{(\psibar A \psi)^{k}}{k+1}
\, (\psibar X B^k Y \psi)
\bigg)
\\
& 
\hspace{4cm}
=
\int \mathcal{D}(\psi, \psibar)
\exp \bigg(-
\psibar X \,
\frac{\ln (1 - (\psibar A \psi) B)}{(\psibar A \psi) B}
\, Y \psi
\bigg)
\ef.
\end{split}
\ee
\item[(b)] If $X$ and $Y$ are column-pseudo-commutative,
then
\be
\begin{split}
&
\hspace{-10mm}
\sum_{\substack{L \subseteq [m] \\ |L|= n}} 
\rowdet X_{[n],L} \rowdet Y_{L,[n]}
=
\int \mathcal{D}(\psi, \psibar)
\exp \bigg(
\sum_{k \geq 0}
\, (\psibar X B^k Y \psi) \, \frac{(\psibar A \psi)^{k}}{k+1}
\bigg)
\\
& 
\hspace{3.5cm}
=
\int \mathcal{D}(\psi, \psibar)
\exp \bigg(-
\psibar X \,
\frac{\ln (1 - (\psibar A \psi) B)}{(\psibar A \psi) B}
\, Y \psi
\bigg)
\ef.
\end{split}
\ee
\end{itemize}
\end{prop}
The commutation condition (\ref{eq.XAgrass}) in the hypotheses above
is identical to the condition (\ref{XA2}) in Proposition
\ref{new.capCB}. Thus, as stated earlier, the following Lemmas
\ref{lem.XABcond} and \ref{lem.XABcond2} discuss mild conditions on
$B$ that would imply it.

However we are not aware of equally satisfactory conditions under
which the hypothesis (\ref{eq.YAgrass}) holds. In particular, the
hypothesis that $Y$ is row-commutative would have rather suggested to
interchange indices $i$ and $k$ in the second summand, instead of $j$
and $\ell$.  A sufficient condition would be that $Y$ is both row- and
column-pseudo-commutative, i.e., that it is \emph{tout-court}
commutative, as in this situation the column-analogue of Lemmas
\ref{lem.XABcond} and \ref{lem.XABcond2} would apply (note, with the
hypotheses of the lemmas now being on $B^{\rm T}$).  We are not aware
of any set of matrices realizing the hypotheses of the proposition
above and in which $Y$ is not commutative, nor we have a proof that
such a realization cannot exist (see the discussion at the end of
Section~\ref{sec.NCCB}).

\vspace{2mm}

\noindent
We will prove Proposition~\ref{prop.old} in Section~\ref{sec.equiv}.
Then in Section~\ref{sec.NCCB} we recall some basic facts which were
useful in our proof of Proposition~\ref{theo.capCB}, and will also be
needed in the following. This section includes also a discussion on
the conditions on the commutation of $X$ and $A$.
Section~\ref{sec.enumeration} is of combinatorial nature. It presents
a lemma on the weighted enumeration of a family of lattice paths, (of
\emph{{\L}ukasiewicz} type), that is used later on in our proofs of
Capelli-like identities.  Section~\ref{sec.full} presents the proof of
Proposition~\ref{new.capCB}, the non-commutative Cauchy--Binet formula
in Quantum oscillator algebra representation.  Section~\ref{sec.holom}
presents a small variant of this formla, in which coherent states of
the quantum oscillator are used.
In Section~\ref{sec.cbh} we derive a useful specialization of the
Campbell-Baker-Hausdorff formula, which we use in
Section~\ref{sec.grass} to give a proof of
Proposition~\ref{prop:grass}, the non-commutative Cauchy--Binet
formula in Grassmann Algebra representation. In Section
\ref{sec.directGrass} we give a short proof of 
Proposition~\ref{prop:grass}, for the case $B=I$.

\section{The bosonic oscillator and multilinear non-commutative
  functions}
\label{sec.equiv}

At the beginning of Section \ref{ssec.bosointro}, we set some
notations for the bosonic oscillator. Among other things, we fixed the
normalization of the state $\ket{0}$. There exists a residual freedom
in choosing the relative norm of states $\ket{n}$, that we fix here,
by setting for each $m, n\in \mathbb{N}$
\begin{align}
\big(\adg\big)^n  
\ket{m}
&= 
\ket{m+n}
\ef,
&
\bra{m}
\aaa^n 
&=
\bra{m+n}
\ef,
\label{adagger}
\end{align}
from which it follows
\begin{align}
\aaa^n 
\ket{m}
&=
\frac{m!}{(m-n)!} \, 
\ket{m-n}
\ef,
&
\bra{m}
\big(\adg\big)^n 
&=
\bra{m-n}
\frac{m!}{(m-n)!}
\ef,
\label{a}
\end{align}
and
\be
\braket{n}{m}
= n!
\, \delta_{nm} 
\label{norma}
\ef.
\ee
As, for $m \in \mathbb{N}$, the states $|m\rangle$ form a complete
set, we have
\be
\label{eq.resoId}
1
=
\sum_{m \geq 0}
\ket{m}
\frac{1}{m!}
\bra{m}
\ef,
\ee
as operators acting on the Hilbert space.

In this section we prove Proposition~\ref{prop.old}.
The two cases are analogous, and we study the `row' case, that is
we choose to prove identity~\reff{old.row}.
We shall in fact prove a more general result, for a family of
multilinear non-commutative functions. Both results are statements on
the fact that, taking scalar products, implement substitutional rules
on suitable polynomials in the algebra of the quantum oscillator,
in a way non dissimilar to the
content of `modern' umbral calculus \emph{a'la} Rota.

\begin{prop}
\label{theo.multilin}
Let $R$ be a ring, $k$, $n$ and $\{m(i)\}_{1 \leq i \leq n}$ integers,
and $\{x_{ij}^{(h)}\}$ a collection of expressions in $R$, for $0 \leq
h \leq k$, $i\in [n]$ and $j\in [m(i)]$.  Consider also a
Weyl-Heisenberg algebra as in~(\ref{WH}), with operators commuting
with the $x$'s.  Take $f(\aaa)$ a formal power series in $\aaa$, such
that $f(0)=1$ and $f'(0) \neq 0$, so that both $f(\aaa)$ and
$f'(\aaa)$ are invertible.
Consider a further indeterminate $s$, and let $g(\aaa,s)$ be the formal
power series in $\aaa$ and $s$ defined as
\be
\label{eq.defginv}
g(\aaa,s) := s\, \left[ \frac{\partial}{\partial \aaa} f(\aaa)^{-s}\right]^{-1}
= - [f'(\aaa)]^{-1} f(\aaa)^{s+1}
\ef.
\ee
Then, introduce the operators
\be
\chi_h (\aaa, \adg)
:=
\frac{1}{h!}
\big( \adg g(\aaa,s) \big)^{h}
f(\aaa)^{-sh-1}
\ef.
\ee
Let
\be
y_{ij} :=  \sum_{h=0}^{k} \binom{i-1}{h}_{\!s} x_{ij}^{(h)}
\ef
\ee
with
\be
\binom{\ell}{h}_{\!s}
: =
\frac{1}{h!} \ell (\ell-s) \cdots (\ell-(h-1)s)
=
\left\{
\begin{array}{ll}
s^h \binom{\ell/s}{h} & s \neq 0 \ef; \\
\frac{\ell^h}{h!}     & s = 0 \ef.
\end{array}
\right.
\ee
Define
\be
z_{ij}(\aaa,\adg) :=
\sum_{h=0}^{k} \chi_h(\aaa,\adg) x_{ij}^{(h)}
\ef.
\ee
Then, for any polynomial $\phi$ of the $\mathcal{N}$ variables
$\{y_{ij}\}$ in the ring $R$, homogeneous of degree $n$, 
and with monomials of the form $\prod_{i=1}^n y_{ij(i)}$ (with the
product in order),\;\footnote{This means that $\phi$ is
  \emph{multilinear} in each set $Y_i = \{ y_{ij} \}_{j\in [m(i)]}$.}
the following representation holds
%
%
\be
\label{eq.qoa}
\phi(\{y_{ij}\})
=
\bra{0} \phi \big(\{z_{ij}(\aaa,\adg) \} \big)  \ket{0}
\ef.
\ee
\end{prop}

\noindent
We recognize the identity~(\ref{old.row}) as a special case, with
$k=1$, $x^{(0)}_{ij} = U_{ij}$, $x^{(1)}_{ij} = V_{ij}$ and $f(\aaa) =
1-\aaa$.  (Thus in particular, $\chi_0=(1-\aaa)^{-1}$ and $\chi_1 =
\adg$.) The polynomial $\phi$ is chosen to be $\phi(y) = \rowdet Y$,
for $Y$ the matrix with entries $y_{ij} = U_{ij} + (i-1) V_{ij} $.
This correspondence is valid regardless of $s$, as $s$ appears
explicitly only for $k \geq 2$.

Towards the end of the proof of this theorem we will need a Lemma
in quantum oscillator algebra, which we prove immediately
\begin{lemma} 
\label{lem.faf}
For any indeterminates $\ell$ and $s$, $f(\aaa)$ and $g(\aaa,s)$ as
above, and any $h$ and $m$ in $\mathbb{N}$,
\be
C_{\ell,h,m} 
:= \frac{1}{h!} \bra{0} 
f(\aaa)^{-\ell}
\big( \adg g(\aaa,s) \big)^h 
f(\aaa)^{\ell-hs} \ket{m}
= \binom{\ell}{h}_{\!s} \delta_{m,0}
\ef.
\ee
\end{lemma}
\proof Indeed, if $h=0$ we trivially have 
$C_{\ell,0,m} = \braket{0}{m} = \delta_{m,0}$, while
if $h>0$ we can write
\be
\begin{split}
C_{\ell,h,m} 
& =
\frac{1}{h!}
\bra{0} 
f(\aaa)^{-\ell}
\adg
g(\aaa,s)
\big( \adg g(\aaa,s) \big)^{h-1}
f(\aaa)^{\ell-hs} \ket{m}
\\
& =
\frac{1}{h!}
\bra{0} 
\Big(
\adg
f(\aaa)^{-\ell} 
+
\big[ f(\aaa)^{-\ell},\adg \big]
\Big)
g(\aaa,s)
\big( \adg g(\aaa,s) \big)^{h-1}
f(\aaa)^{\ell-hs} \ket{m}
\\ & =
\frac{\ell}{h!}
\bra{0} 
f(\aaa)^{-\ell-1}
(-f'(\aaa) g(\aaa,s))
\big( \adg g(\aaa,s) \big)^{h-1}
f(\aaa)^{\ell-hs} \ket{m}
\\ & =
\frac{\ell}{h!}
\bra{0} 
f(\aaa)^{-(\ell-s)}
\big( \adg g(\aaa,s) \big)^{h-1}
f(\aaa)^{(\ell-s)-(h-1)s} \ket{m}
\\ & = \frac{\ell}{h} C_{\ell-s,h-1,m}
\ef,
\end{split}
\ee
where we used the fact that $\bra{0} \adg = 0$, and the definition
(\ref{eq.defginv}). So we get the result by induction in $h$.
\hfill \qed


\proofof{Proposition~\ref{theo.multilin}}
A generic monomial of $\phi$ can be
labeled by a vector $J=(j(1), \ldots, j(n)) \in [m(1)] \times \cdots
\times [m(n)]$, thus $\phi$
has the form 
\be
\phi(\{y_{ik}\}) = \sum_{J= (j(1), \ldots, j(n))} 
c_J \; \prod_{i=1}^n y_{i j(i)} 
\ef.
\ee
Both $y_{ij}$'s and $z_{ij}(\aaa,\adg)$'s are defined as a sum of
$k+1$ terms. Perform the corresponding expansion on both sides of
(\ref{eq.qoa}), and label each term by a vector $\mu \in \{ 0, \ldots,
k \}^n$. For the expression on the left hand side we have
\be
\phi(\{y_{ik}\}) = \sum_{J, \, \mu} c_J
\bigg(
\prod_{i=1}^n
\binom{i-1}{\mu(i)}_{\!s}
\bigg)
\prod_{i=1,\ldots,n} 
x^{(\mu(i))}_{i j(i)}
\ef,
\ee
while for the one on the right hand side we have
\be
\bra{0}
\phi
\big(
\big\{
z_{ik}(\aaa,\adg)
\big\}
\big)
\ket{0}
=
\sum_{J,\,\mu} c_J
\bra{0}
\!\!
\prod_{i=1,\ldots,n}
\!\!
\chi_{\mu(i)}
\ket{0}
\prod_{i=1,\ldots,n}
\!\!
x_{i j(i)}^{(\mu(i))}
\ef.
\ee
As the $x_{ij}^{(h)}$ are arbitrary non-commuting indeterminates, and
$\phi$ is arbitrary, the
identity must hold separately for each summand labeled by a pair
$(J,\mu)$, i.e.~that for any vector~$\mu$ we have to prove that
\be
\prod_{\substack{\ell\in [n] \\ \mu(\ell)\neq 0}}
\binom{\ell-1}{\mu(\ell)}_{\!s}
=
\bra{0}
\!\!
\prod_{i=1,\ldots,n}
\!\!
\chi_{\mu(i)}
\ket{0}
\ef.
\ee
Let $(\ell_1, \ldots, \ell_k)$ be the ordered list of indices
$i$ such that $\mu(i) \neq 0$, so that
\be
\prod_{i=1,\ldots,n}
\!\!
\chi_{\mu(i)}
=
\chi_0^{\ell_1-1} \, \chi_{\mu(\ell_1)} \,
\chi_0^{\ell_2-\ell_1-1} \, \chi_{\mu(\ell_2)} \,
\chi_0^{\ell_3-\ell_2-1} \, \chi_{\mu(\ell_3)} \,
\cdots
\, \chi_{\mu(\ell_k)} \, \chi_0^{n-\ell_k} 
\ee
where all the powers are non-negative integers, and all
$\mu(\ell_j)$'s are in the range $\{1,\ldots,k\}$.
The expression $\chi_0^{-1} = f(a)$ is defined as a formal power
series, and we can write
\be
\prod_{i=1,\ldots,n}
\!\!
\chi_{\mu(i)}
=
\Big(
\prod_{\alpha=1,\ldots,k}
\!\!
\chi_0^{\ell_{\alpha}-1} \, \chi_{\mu(\ell_{\alpha})} \, \chi_0^{-\ell_{\alpha}} \;
\Big)
\chi_0^{n}
\ef.
\ee
Let
\be
\hat{O}_{\ell, h}
:=
\chi_0^{\ell-1}
\chi_h
\chi_0^{-\ell}
\ef,
\ee
we need to prove that, for any $k$-uple $\ell_1 < \cdots < \ell_k$,
\be
\prod_{\alpha=1}^{k}
\binom{\ell_{\alpha} - 1}{\mu(\ell_{\alpha})}_{\!s} =
\bra{0} 
\bigg(
\prod_{\alpha=1,\ldots,k} \!\!
\hat{O}_{\ell_{\alpha}, \mu(\ell_{\alpha}) }
\bigg)
f(\aaa)^{-n}
\ket{0}
\ef.
\ee
First of all realize that $f(\aaa)^{-n} \ket{0} = \ket{0}$. Then,
because of Lemma~\ref{lem.faf},
\be
\bra{0}  \hat{O}_{\ell_{1}, \mu(\ell_{1}) } \ket{m} = \delta_{m,0}
\binom{\ell_{1} - 1}{\mu(\ell_{1})}_{\!s}
\ee
so that by introducing a resolution of the identity,
equation (\ref{eq.resoId}), we get a recursion in~$\alpha$
\be
\begin{split}
\bra{0} 
\prod_{\alpha=1,\ldots,k} \!\!
\hat{O}_{\ell_{\alpha}, \mu(\ell_{\alpha}) }
\ket{0}
&=
\sum_{m \geq 0}
\bra{0} 
\hat{O}_{\ell_{1}, \mu(\ell_{1}) }
\ket{m} \frac{1}{m!}
\bra{m}
\prod_{\alpha=2,\ldots,k} \!\!
\hat{O}_{\ell_{\alpha}, \mu(\ell_{\alpha}) }
\ket{0}
\\
&  =
\sum_{m \geq 0}
\frac{\delta_{m,0}}{m!}
\binom{\ell_{1} - 1}{\mu(\ell_{1})}_{\!s}
\bra{m}
\prod_{\alpha=2,\ldots,k} \!\!
\hat{O}_{\ell_{\alpha}, \mu(\ell_{\alpha}) }
\ket{0}
\\
&  =
\binom{\ell_{1} - 1}{\mu(\ell_{1})}_{\!s}
\bra{0}
\prod_{\alpha=2,\ldots,k} \!\!
\hat{O}_{\ell_{\alpha}, \mu(\ell_{\alpha}) }
\ket{0}
\ef,
\end{split}
\ee
which proves the statement of the theorem.
\qed


\section{Some properties of commutators}
\label{sec.NCCB}


Let us begin by recalling two elementary 
facts~\cite[Lemma 2.1 and 2.2]{uscapelli} that we used repeatedly and
shall use in this paper:

\begin{lemma}[Translation Lemma]
  \label{lemma.translation}
Let $\scra$ be an abelian group, and let $f \colon\, \scrs_n \to \scra$.
Then, for any $\tau \in \scrs_n$, we have
\be
   \sum_{\sigma \in \scrs_n}  \sgn(\sigma) \, f(\sigma)  =
   \sgn(\tau) \sum_{\sigma \in \scrs_n}  \sgn(\sigma) \, f(\sigma \circ \tau)
   \;.
\ee
\end{lemma}

\proof
Just note that both sides equal
$\sum\limits_{\sigma \in \scrs_n}
    \sgn(\sigma \circ \tau) \, f(\sigma \circ \tau)$.
\qed

\begin{lemma}[Involution Lemma]
  \label{lemma.involution}
Let $\scra$ be an abelian group, and let $f \colon\, \scrs_n \to \scra$.
Suppose that there exists a pair of distinct elements $i,j \in [n]$
such that
\be
     f(\sigma)  =  f(\sigma \circ (ij))
\ee
for all $\sigma \in \scrs_n$
[where $(ij)$ denotes the transposition interchanging $i$ with~$j$].
Then
\be
   \sum_{\sigma \in \scrs_n}  \sgn(\sigma) \, f(\sigma) = 0
   \;.
\ee
\end{lemma}

\proof
We have
\begin{align}
   \sum_{\sigma \in \scrs_n}  \sgn(\sigma) \, f(\sigma)
    = &
   \sum_{\sigma \colon\, \sigma(i) < \sigma(j)}  \! \sgn(\sigma) \, f(\sigma)
   \;+\!\!
   \sum_{\sigma \colon\, \sigma(i) > \sigma(j)}  \! \sgn(\sigma) \, f(\sigma)
       \\[1mm] 
    = &
   \sum_{\sigma \colon\, \sigma(i) < \sigma(j)}  \! \sgn(\sigma) \, f(\sigma)
   \;-\!\!
   \sum_{\sigma' \colon\, \sigma'(i) < \sigma'(j)}  \!\! \sgn(\sigma') \,
     f(\sigma' \circ (ij))
       \qquad \\[1mm] 
    = &
   \quad 0  \;,
\end{align}
where in the second line we made the change of variables
$\sigma' = \sigma \circ (ij)$ and used $\sgn(\sigma') = -\sgn(\sigma)$
[or equivalently used the Translation Lemma].
\qed

In the following we shall need of a less restrictive notion than the
pseudo-commu\-tative matrix. Let us begin by observing that $\mu_{ijkl}
\bydef [M_{ij}, M_{kl}]$ is manifestly antisymmetric under the
simultaneous interchange $i \leftrightarrow k$, $j \leftrightarrow l$.
So symmetry under one of these interchanges is equivalent to
antisymmetry under the other.  Let us therefore say that a matrix $M$
has {\em row-symmetric}\/ (and {\em column-antisymmetric}\/) {\em
  commutators}\/ if $[M_{ij}, M_{kl}] = [M_{kj}, M_{il}]$ for all
$i,j,k,l$, and {\em column-symmetric}\/ (and {\em
  row-antisymmetric}\/) {\em commutators}\/ if $[M_{ij}, M_{kl}] =
[M_{il}, M_{kj}]$ for all $i,j,k,l$.

Then we shall need the following two lemmas.
\begin{lemma} 
\label{lem.XXcomm}
For a $n$-dimensional matrix $M$ with row-symmetric  commutators, that is
satisfying
\be
\label{eq.XX868565}
[M_{ij}, M_{kl} ] - [M_{kj}, M_{il} ] = 
0 \qquad 
\hbox{for all\ }\quad  i,\, j,\, k,\, l
\ef,
\ee
any vector $(\ell_1, \ldots, \ell_n)$,
and any permutation $\pi \in \mathcal{S}_n$,
\be
\sum_{\sigma \in \scrs_n}
\sgn(\sigma) \prod_{i=1}^n
M_{\sigma(i)\, \ell_{i}}
=
\sum_{\sigma \in \scrs_n}
\sgn(\sigma) \prod_{i=1}^n
M_{\sigma \pi(i)\, \ell_{\pi(i)}}
\ef.
\ee
\end{lemma}
It suffices to prove the lemma for a single transposition of elements, 
consecutive after the permutation $\sigma$, namely
$\pi = (\sigma(i)\, \sigma(i+1))$. We denote as $L_\sigma$ and
$R_\sigma$ the factors on left and on the right (note that they do not depend
from $\sigma(i)$ and $\sigma(i+1)$).  We can write the statement as
\be
\sum_{\sigma \in \scrs_n}
\sgn(\sigma)
L_\sigma
M_{\sigma(i)\,   \ell_{i}}
M_{\sigma(i+1)\, \ell_{i+1}}
R_\sigma
=
\sum_{\sigma \in \scrs_n}
\sgn(\sigma)
L_\sigma
M_{\sigma(i+1)\, \ell_{i+1}}
M_{\sigma(i)\,   \ell_{i}}
R_\sigma
\ef.
\ee
The difference of the two expressions is, by definition,
\be
\sum_{\sigma \in \scrs_n}
\sgn(\sigma)
L_\sigma
[M_{\sigma(i)\,   \ell_{i}} ,
M_{\sigma(i+1)\, \ell_{i+1}} ]
R_\sigma
\ee
which vanishes because the hypothesis~\reff{eq.XX868565} allows the
application of the Involution Lemma. \qed

Now we have a sequence of lemmas exploring the consequences of the
Jacobi identity.
\begin{lemma}
\label{lemma:B}
Let $R$ be a 
ring,
and let $X$ 
and $Y$ 
be matrices with elements in $R$. 
\begin{itemize}
\item[(a)] If  $X$ is row-pseudo-commutative
then, at fixed $Y_{ef}$, for all $a,b,c,d$ the antisymmetric part of
$[X_{ab},[X_{cd},Y_{ef}]]$ in the exchange of $a$  with $c$ is
symmetric in the exchange of $b$ and $d$, that is
\begin{gather}
[X_{ab},[X_{cb},Y_{ef}]] - [X_{cb},[X_{ab},Y_{ef}]]
= 0
\label{XA1pre}
\\
[X_{ab},[X_{cd},Y_{ef}]] - [X_{cb},[X_{ad},Y_{ef}]]
+ [X_{ad},[X_{cb}, Y_{ef}]] - [X_{cd}, [X_{ab}, Y_{ef}]] 
= 0 
\label{XA1}
\ef.
\end{gather}
\item[(b)] If  $Y$ is column-pseudo-commutative then, at fixed
  $X_{ef}$, for all $a,b,c,d$ the antisymmetric part of
  $[Y_{ab},[Y_{cd},X_{ef}]]$ in the exchange of $b$  with $d$ is
  symmetric in the exchange of $a$ and $c$, that is
\begin{gather}
[Y_{ab},[Y_{ad},X_{ef}]] -  [Y_{ad},[Y_{ab},X_{ef}]] 
=
0
\\
[Y_{ab},[Y_{cd},X_{ef}]] -  [Y_{ad},[Y_{cb},X_{ef}]]   +  [Y_{cb},
  [Y_{ad}, X_{ef}]] -  [Y_{cd}, [Y_{ab}, X_{ef}]] = 0 
\ef.
\end{gather}
%
\end{itemize}
\end{lemma}
\proof
(a) Start from the Jacobi Identity applied to the triplet $(X_{ab}, X_{cd}, Y_{ef})$,
\be
[X_{ab},[X_{cd},Y_{ef}]] + [Y_{ef}, [X_{ab},X_{cd}]]  + [X_{cd}, [Y_{ef}, X_{ab}]] = 0
\ee
If we set $d=b$, as $X$ is row-pseudo-commutative, $[X_{ab},X_{cb}]=0$
so that (\ref{XA1pre}) follows.
For (\ref{XA1}), consider also the Jacobi identity for the triplet
$(X_{cb}, X_{ad}, Y_{ef})$ to obtain
\be
  [X_{cb}, [X_{ad}, Y_{ef}]] 
+ [Y_{ef}, [X_{cb}, X_{ad}]]
+ [X_{ad}, [Y_{ef}, X_{cb}]] = 0
\ee
so that, by subtraction and the hypothesis that $X$ is
row-pseudo-commutative then
\be
[X_{ab},[X_{cd},Y_{ef}]] -  [X_{cb},[X_{ad},Y_{ef}]]   +  [X_{ad},
  [X_{cb}, Y_{ef}]] -  [X_{cd}, [X_{ab}, Y_{ef}]] = 0
\ef.
\ee
The proof of (b) is similar. \qed

This lemma implies the following
\begin{coroll}
\label{coro.XABpre}
If $X$ and $Y$ are as in case (a) of the lemma above, and furthermore
they satisfy the commutation relation (\ref{new}),
$[X_{ij},Y_{k\ell }] =  \, -A_{i\ell } B_{kj}$, then, for every
$a, b, c, e, f$,
\be
\label{eq.incor1}
([X_{ab},A_{cf}] - [X_{cb},A_{af}]) B_{eb}
=0
\ef;
\ee
and for every $a, b, c, d, e, f$,
\be
\label{eq.incor2}
([X_{ab},A_{cf}] - [X_{cb},A_{af}]) B_{ed}
+
(b \leftrightarrow d)
=
0
\ef.
\ee
\end{coroll}

\noindent
We are now ready to state sufficient conditions on $B$, for having
the commutation relation (\ref{XA2}), 
$[X_{ij}, A_{k\ell}] - [X_{kj}, A_{i\ell}] = 0$.

Recall that, in a ring $R$, a nonzero element $x$ is a \emph{left zero
  divisor} if there exists a nonzero $y$ such that $xy=0$.  Right zero
divisors are analogously defined.  A nonzero element of a ring that is
not a left zero divisor is called \emph{left-regular} (and
analogously for right). Then
\begin{lemma}
\label{lem.XABcond}
Let $X$, $A$ and $B$ as in Corollary \ref{coro.XABpre}, of sizes
respectively $n \times m$, $n \times n$ and $m \times m$, and $B_{ij}$
commuting with every other matrix element. Suppose that there exist an
index $d \in [m]$, and a vector $\vec{u} \in R^m$, such that 
$(\vec{u} B)_d$ is
left-regular.
Then
\be
[X_{ij}, A_{k\ell}] - [X_{kj}, A_{i\ell}] = 0
\qquad
\textrm{for all $i,j,k,\ell$.}
\ee
\end{lemma}
\proof
Equations (\ref{eq.incor1}) and (\ref{eq.incor2}) are valid with $B$
written on the left or on the right, as it commutes with
everything.
Consider equation (\ref{eq.incor1}), with arbitrary $a, c, f$, setting
$b=d$, and summing over $e$, after multiplying on the left by $u_e$.
This gives
\be
\bigg( \sum_e u_e B_{ed} \bigg)
([X_{ad},A_{cf}] - [X_{cd},A_{af}]) 
=0
\ef.
\ee
As $(\vec{u} B)_d$ is
left-regular, we obtain 
$[X_{ad},A_{cf}] - [X_{cd},A_{af}] =0$.
Now consider any other index $b \neq d$, and equation
(\ref{eq.incor2}), again summing over $e$, after multiplying on the
left by $u_e$. We obtain
\be
\bigg( \sum_e u_e B_{ed} \bigg)
([X_{ab},A_{cf}] - [X_{cb},A_{af}])
=
-
\bigg( \sum_e u_e B_{eb} \bigg)
([X_{ad},A_{cf}] - [X_{cd},A_{af}])
\ef.
\ee
As the right-most factor on the right hand side is zero, the whole 
right hand side vanishes. As the left-most factor on the left hand
side is left-regular, we have that
$[X_{ab},A_{cf}] - [X_{cb},A_{af}] =0$,
thus completing the proof.
\qed

Furthermore, we can also state
\begin{lemma}
\label{lem.XABcond2}
Let $X$, $A$ and $B$ as in Corollary \ref{coro.XABpre}, of sizes
respectively $n \times m$, $n \times n$ and $m \times m$, and $B_{ij}$
commuting with every other matrix element. Suppose that there exist a
vector $\vec{u} \in R^m$, and a vector $\vec{v} \in R^m$, with $v_i$'s
commuting with $X$, $A$ and $B$ elements and among themselves,
such that the scalar product $2 (\vec{u}, B \vec{v})$ is
left-regular.
Then
\be
[X_{ij}, A_{k\ell}] - [X_{kj}, A_{i\ell}] = 0
\qquad
\textrm{for all $i,j,k,\ell$.}
\ee
\end{lemma}
\proof
Remark that, except for the annoying factor $2$, this lemma is a
generalization of Lemma \ref{lem.XABcond}, to which it (almost) reduces for
$\vec{v}_i = \delta_{i,d}$.

Analogously to Lemma \ref{lem.XABcond},
consider equation (\ref{eq.incor2}), with arbitrary $a, c, f$, summing
over $e, b, d$, after multiplying on the left by $u_e v_b v_d$.
This gives
\be
\begin{split}
&
\bigg( \sum_{e,d} u_e B_{ed} v_d \bigg)
\sum_b ([X_{ab} v_b,A_{cf}] - [X_{cb} v_b,A_{af}]) 
\\
&
\quad
+
\bigg( \sum_{e,b} u_e B_{eb} v_b \bigg)
\sum_d ([X_{ad} v_d,A_{cf}] - [X_{cd} v_d,A_{af}]) 
=0
\ef.
\end{split}
\ee
Performing the sums shows that the two terms are identical.
As $2 (\vec{u}, B \vec{v})$ is
left-regular, we obtain 
$[(X \vec{v})_{a},A_{cf}] - [(X \vec{v})_{c},A_{af}] =0$.
Now take any index $b$, and consider again equation
(\ref{eq.incor2}), but summing only over $e$ and $d$, after multiplying on the
left by $u_e v_d$. We obtain
\be
\bigg( \sum_e u_e B_{ed} v_d \bigg)
([X_{ab},A_{cf}] - [X_{cb},A_{af}])
=
-
\bigg( \sum_e u_e B_{eb} \bigg)
([(X \vec{v})_{a},A_{cf}] - [(X \vec{v})_{c},A_{af}])
\ef.
\ee
As the right-most factor on the right hand side is zero, the whole 
right hand side vanishes. As the left-most factor on the left hand
side is left-regular, we have that
$[X_{ab},A_{cf}] - [X_{cb},A_{af}] =0$,
thus completing the proof.
\qed

An analysis similar to the one of Corollary \ref{coro.XABpre},
performed on matrix $Y$ assumed to be row-pseudo-commutative (remark
that Lemma \ref{lemma:B}(a) exchanging $X$ and $Y$ is a valid starting
point at this aim), gives
\begin{gather}
[Y_{ab},A_{eb}] B_{cf}
=
[Y_{cb},A_{eb}] B_{af}
\ef;
\\
[Y_{ab},A_{ed}] B_{cf} + (b \leftrightarrow d)
=
[Y_{cb},A_{ed}] B_{af} + (b \leftrightarrow d)
\ef.
\end{gather}
These equations are comparatively weaker w.r.t.\ equations
(\ref{eq.incor1}) and (\ref{eq.incor2}), at the aim of establishing
sufficient conditions on $B$ for the hypothesis (\ref{eq.YAgrass}) in
Proposition \ref{prop:grass} to hold. Indeed, while in the previous case we
have already the appropriate exchange structure, mixed to further
exchanges, in this new case the exchange of indices has nothing in
common with (\ref{eq.YAgrass}).

A simple sufficient condition is that $Y$ is in fact commutative,
$[Y_{ij},Y_{k\ell}]=0$ for all $i,j,k,\ell$, as this would imply in
particular that it is column-pseudo-commutative, and the validity
follows from the cases (b) of the lemmas above. Another case leading
to interesting simplifications is when $B$ is the identity matrix, and
$m \geq 2$. In this case, taking $f=c\neq a$
gives
\begin{gather}
[Y_{ab},A_{eb}]
=
0
\ef;
\\
[Y_{ab},A_{ed}] + (b \leftrightarrow d)
=
0
\ef.
\end{gather}
Thus we see that, in this case, either the field has characteristic 2,
or the only possibility for (\ref{eq.YAgrass}) to hold is that
$[Y_{ij},A_{k\ell}]=0$ for all $i$, $k$ and $j \neq \ell$.

\section{A weighted enumeration of {\L}ukasiewicz paths}
\label{sec.enumeration}

Let $n$ an integer. For $0 \leq t \leq n$, consider the
`symbols'
$\vec{\nu}_{t} = \nunu{\nu_1, \cdots, \nu_{t} | 
\nu_{t+1}, \cdots, \nu_n }$,
$n$-uples of integers
with $\nu_i \geq -1$ for $i \leq t$ and 
$\nu_i \geq 0$ for $i > t$.
These symbols are intended as formal indeterminates generating a linear space
over $\mathbb{Z}$.
Consider the quotient given by the relations
\be
\label{eq.symbrule}
\begin{split}
&
\nunu{\nu_1, \cdots, \nu_{t-1} | \nu_t, \cdots, \nu_n }
=
\nunu{\nu_1, \cdots, \nu_{t-1}, \nu_t | \nu_{t+1}, \cdots, \nu_n }
\\
& \hspace{3.5cm}
+
\sum_{k=t+1}^n
\nunu{\nu_1, \cdots, \nu_{t-1},
-1
| \nu_{t+1}, \cdots,
\underbrace{\nu_k+\nu_t+1}_{k\textrm{-th}}, \cdots, \nu_n }
\ef.
\end{split}
\ee
Remark that the sum $|\vec{\nu}_t| = \nu_1 + \cdots + \nu_n$, that we
call the \emph{norm} of the symbol, is homogeneous in all the terms of
the relation, and that, if the left hand side of (\ref{eq.symbrule})
satisfies the bounds above on the $\nu_i$'s, the bounds are satisfied
also by all the summands on the right hand side.

Let us call \emph{height} of 
$\nunu{\nu_1, \cdots, \nu_t | \nu_{t+1}, \cdots, \nu_n }$ the integer
$H=\nu_{t+1} + \cdots + \nu_n$. Then the other combination $\nu_1 +
\cdots + \nu_{t}$ is just the norm minus the height.
We shall call $t$ the {\em level} of
$\nunu{\nu_1, \cdots, \nu_t | \nu_{t+1}, \cdots, \nu_n }$.
We define $V_{t,s}$ as the space of all symbols with level $t$ and
norm~$s$.

Consider any triplet $(t,t',s)$ with $0 \leq t \leq t' \leq n$ and $s
\geq -t$. The relation (\ref{eq.symbrule}) can be seen as a recursion,
allowing to write any symbol $\vec{\nu}_{t} \in V_{t,s}$ as a linear
combination of symbols $\vec{\nu}'_{t'} \in V_{t',s}$.
We will restrict our attention to the symbols with zero norm.
For $t=0$, we have a unique possible symbol in $V_{0,0}$, that is,
$\vec{\nu}_{\bm 0} = \nunu{\, | 0 \cdots 0 }$.
As a consequence, and from the closure property above,
for each $0 \leq t \leq n$
there exists a set of 
integers $c(\vec{\nu}_{t})$
such that
\be
\label{eq.nuZero}
\vec{\nu}_{\bm 0} =
\sum_{\vec{\nu}_{t} \in V_{t,0}}
c(\vec{\nu}_{t}) 
\;
\vec{\nu}_{t}
\ef.
\ee
In the following Lemma \ref{main} we determine a formula for
$c(\vec{\nu}_{t})$, which is the main result of the section.
Before going to the lemma, it is useful to introduce a graphical
interpretation for these symbols.


Symbols of maximal level, 
$\vec{\nu}_n= \nunu{\nu_1, \cdots \nu_n | \,}$,
are in bijection with
paths $\gamma$ on the half-line, that is, if represented as a
`time trajectory' in two dimensions, paths
with height remaining always non-negative, starting at $(0,0)$ and
arriving at $(n,0)$, and with steps of the form $(1,s)$.
The bijection just consists in performing a jump of $-\nu_i$ at
the $i$-th step.  Thus, in our problem we have only steps $s \leq 1$. 
Paths with exactly this set of allowed steps are known as
{\em {\L}ukasiewicz paths} (see~\cite[pag.~71]{Flajolet-book} or
\cite[Example 3, pag.~14]{Cyril}).  An example of symbol-path
correspondence is
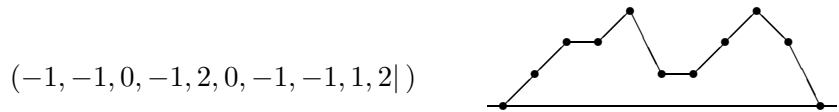
\[
\nunu{-1,-1,0,-1,2,0,-1,-1,1,2|\,}
\hspace{1cm}
\setlength{\unitlength}{12pt}
\raisebox{-3mm}{
\begin{picture}(10,3.4)(0,-0.1)
\put(-0.5,0){\line(1,0){11}}
\put(0,0){\line(1,1){1}}
\put(1,1){\line(1,1){1}}
\put(2,2){\line(1,0){1}}
\put(3,2){\line(1,1){1}}
\put(4,3){\line(1,-2){1}}
\put(5,1){\line(1,0){1}}
\put(6,1){\line(1,1){1}}
\put(7,2){\line(1,1){1}}
\put(8,3){\line(1,-1){1}}
\put(9,2){\line(1,-2){1}}
\put(0,0){\circle*{0.25}} 
\put(1,1){\circle*{0.25}} 
\put(2,2){\circle*{0.25}} 
\put(3,2){\circle*{0.25}} 
\put(4,3){\circle*{0.25}} 
\put(5,1){\circle*{0.25}} 
\put(6,1){\circle*{0.25}} 
\put(7,2){\circle*{0.25}} 
\put(8,3){\circle*{0.25}} 
\put(9,2){\circle*{0.25}} 
\put(10,0){\circle*{0.25}}
\end{picture}
}
\]
More generally, symbols of level $t$ and height $H$ are in bijection
with pairs $(\gamma,\pi)$, where $\gamma$ is a path as above,
terminating at $(t,H)$, and $\pi$ is a partition of $H$ `stones' into
$n-t$ boxes (that we represent graphically as the columns with indices
from $t+1$ to $n$, following the path). For example
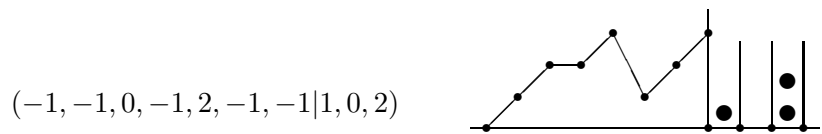
\[
\nunu{-1,-1,0,-1,2,-1,-1|1,0,2}
\hspace{1cm}
\setlength{\unitlength}{12pt}
\raisebox{-3mm}{
\begin{picture}(10,4)(0,-0.2)
\put(-0.5,0){\line(1,0){11}}
\put(0,0){\line(1,1){1}}
\put(1,1){\line(1,1){1}}
\put(2,2){\line(1,0){1}}
\put(3,2){\line(1,1){1}}
\put(4,3){\line(1,-2){1}}
\put(5,1){\line(1,1){1}}
\put(6,2){\line(1,1){1}}
\put(0,0){\circle*{0.25}} 
\put(1,1){\circle*{0.25}} 
\put(2,2){\circle*{0.25}} 
\put(3,2){\circle*{0.25}} 
\put(4,3){\circle*{0.25}} 
\put(5,1){\circle*{0.25}} 
\put(6,2){\circle*{0.25}} 
\put(7,3){\circle*{0.25}} 
\put(7,0){\circle*{0.25}} 
\put(8,0){\circle*{0.25}} 
\put(9,0){\circle*{0.25}} 
\put(10,0){\circle*{0.25}}
\put(7,0){\line(0,1){3.75}}
\put(8,0){\line(0,1){2.75}}
\put(9,0){\line(0,1){2.75}}
\put(10,0){\line(0,1){2.75}}
\put(7.5,0.5){\circle*{0.5}}
\put(9.5,0.5){\circle*{0.5}}
\put(9.5,1.5){\circle*{0.5}}
\end{picture}
}
\]
Paths in one dimension can be described equivalently, either by the
sequence of jumps $-\nu_i$, as above, or by the height profile
$h_i = \sum_{j=1}^{i} (-\nu_j)$.
Both notations will be useful in the following.

One easily sees that a necessary condition for $c(\vec{\nu}_{t}) \neq
0$ is that the corresponding path never goes below the horizontal
axis. Indeed, the recursion is such that, if the left hand side of
(\ref{eq.symbrule}) has non-negative height $H$, then this is true
also for all the summands on the right hand side. Another way of
seeing this property is to realize that our graphical structures
$(\gamma,\pi)$ form a family which is stable under the recursion, and
$H$, which is both the final height in the path and the number of
stones, must remain always non-negative.

Our lemma states
\begin{lemma}
\label{main}
For $\vec{\nu}_{t} = (\gamma,\pi)$,
the function $c(\vec{\nu}_{t})$ 
depends only on $\gamma$ (and not on $\pi$), and
is given by
\be
\label{eq.cgam2}
c(\gamma) = 
h_{t}!
\prod_{\substack{i \in [t] \\ h_i \leq h_{i-1} }} 
\frac{h_{i-1}!}{h_i!}
\ef.
\ee
In particular, when $t = n$, the path must have $h_n=0$ and therefore
\be
\label{eq.cgam}
c(\gamma) = 
\prod_{\substack{i \in [n] \\ h_i \leq h_{i-1} }} 
\frac{h_{i-1}!}{h_i!}
\ef.
\ee
\end{lemma}

\proof
Consider equation (\ref{eq.symbrule}) to derive a recursion for the
coefficients. For the symbol
$\vec{\nu}_t=\nunu{\nu_1, \cdots, \nu_t | \nu_{t+1}, \cdots, \nu_n }$
we have
\begin{equation}
\label{eq.symbruleb}
c\big( \vec{\nu}_t \big)
=
\begin{cases}
c\big(
\nunu{\nu_1, \!\cdots\!, \nu_{t-1} | \nu_t, \!\cdots\!, \nu_n }
\big) & \hbox{if\,}\ \nu_t \geq 0 \ef;\\
\; \rule{0pt}{8mm}%
\displaystyle{
\!\!
\sum_{k=t+1}^n
\sum_{\nu'=1}^{\nu_k}
c\big(
\nunu{
\nu_1, \!\cdots\!, \nu_{t-1} | 
\nu'-1, \nu_{t+1}, \!\cdots\!,
\underbrace{\nu_k - \nu'}_{k\textrm{-th}}, \!\cdots\!, \nu_n }
\big) }
& \hbox{if\,}\ \nu_t = -1
\ef.
\end{cases}
\ee
%
%
We proceed by induction in $t$, starting from the trivial unique
solution $c(\vec{\nu}_{\bm 0})$ of (\ref{eq.nuZero}) for $t=0$.
Assuming the formula for $c(\vec{\nu}_{t})$ valid up to $t-1$, we have
\be
c\big( \vec{\nu}_t \big)
=
\begin{cases}
\displaystyle{h_{t-1}!
\prod_{\substack{i \in [t-1] \\ h_i \leq h_{i-1} }} 
\frac{h_{i-1}!}{h_i!}} & \hbox{if\,}\  \nu_t \geq 0 \ef;
\\
\displaystyle{h_{t-1}!
\prod_{\substack{i \in [t-1] \\ h_i \leq h_{i-1} } }
\frac{h_{i-1}!}{h_i!} \sum_{k=t+1}^n
\sum_{\nu'=1}^{\nu_k}
1 }
& \hbox{if\,}\  \nu_t = -1
\ef.
\end{cases}
\ee
In the case $\nu_t \geq 0$, we have $h_t \leq h_{t-1}$ and therefore
\be
h_{t-1}! \prod_{\substack{i \in [t-1] \\ h_i \leq h_{i-1} } }
\frac{h_{i-1}!}{h_i!} 
= 
h_t! \prod_{\substack{i \in [t] \\ h_i \leq h_{i-1} } } 
\frac{h_{i-1}!}{h_i!}
\ee
as required.
If $\nu_t =-1$, remark that
\be
\sum_{k=t+1}^n \sum_{\nu'=1}^{\nu_k} 1 = \sum_{k=t+1}^n \nu_k = h_t
\ef,
\ee
then, as $h_t = h_{t-1}+1 > h_{t-1}$, we soon get that
\be
h_t \, h_{t-1}! 
\prod_{\substack{i \in [t-1] \\ h_i \leq h_{i-1} } }
\frac{h_{i-1}!}{h_i!} 
= 
h_{t}! \prod_{\substack{i \in [t-1] \\ h_i \leq h_{i-1} } } 
\frac{h_{i-1}!}{h_i!}
= h_{t}! \prod_{\substack{i \in [t] \\ h_i \leq h_{i-1} } } 
\frac{h_{i-1}!}{h_i!}
\ef,
\ee
which completes the proof.\qed

\noindent
Now, for symbols of maximal level, 
$\nunu{\nu_1, \cdots \nu_n | \,}$,
we give a representation in quantum
oscillator algebra of the combinatorial formula for the coefficients 
$c(\vec{\nu}_{t})$

\begin{lemma}
\label{lemma:WH}
For $\nu \geq -1$, define
the operator in the Weyl-Heisenberg algebra
\be
\chi(\nu)= \begin{cases}
(\adg)^{\nu} & 
\nu \geq 0 \ef;
\\
\aaa & 
\nu = -1 \ef.
\end{cases}
\ee
Then, when the symbol $\vec{\nu}_n=\nunu{\nu_1, \cdots \nu_n | \,}$
corresponds to a path $\gamma$ as described above,
\be
\bra{0} \chi(\nu_1) \cdots \chi(\nu_n) \ket{0} = 
c(\vec{\nu}_{n})
=
\prod_{\substack{i \in [n] \\ h_i \leq h_{i-1} }} 
\frac{h_{i-1}!}{h_i!}
\ef,
\ee
while otherwise
\be
\bra{0} \chi(\nu_1) \cdots \chi(\nu_n) \ket{0} = 0
\ef.
\ee 
\end{lemma}

\proof We proceed by induction.  Assume that, for a sequence
$\nu_1,\dots, \nu_t$ such that the corresponding path remains
positive,
\be
\bra{0} \chi(\nu_1) \cdots \chi(\nu_t) = 
\bra{h_t}
\! \prod_{\substack{i \in [t] \\ h_i \leq h_{i-1} }} 
\!\! \frac{h_{i-1}!}{h_i!}
\ef.
\ee
Then, we analyse the application of the operator $\chi(\nu_{t+1})$ to
the right. If $\nu_{t+1} = -1$,
because of~\reff{adagger}, the application
of $\aaa$ consistently brings $\bra{h_t}$ to 
$\bra{h_t + 1} = \bra{h_{t+1}}$.
If $\nu_{t+1} \geq 0$,
because of \reff{a}, the application of $(\adg)^{\nu}$ brings
$\bra{h_t}$ to $\bra{h_t-\nu}$,
with an extra factor $h_t!/(h_t - \nu)!$ (which, in
particular, is zero if the path goes below the horizontal axis).
Taking finally the scalar product with $\ket{0}$ ensures that the path
ends at height zero. \qed

\section{The Capelli identity in Weyl--Heisenberg Algebra}
\label{sec.full}

We are now ready for the:

\proofof{Proposition~\ref{new.capCB}}
(a) As a first step, by simply using the fact that $X$ is
row-pseudo-commutative, in~\cite[Section 3]{uscapelli} we get that
\be
\sum_{\substack{ L \subseteq [m] \\ |L|=n }} 
\coldet X_{[n],L} \coldet Y_{L,[n]} 
=
\sum_{\sigma \in \scrs_n} \sgn(\sigma)
\!
\sum_{l_{1},\dots, l_{n}\in [m]}
\!
\bigg( \prod_{i=1}^n X_{\sigma(i)\, l_{i}} \bigg)
\prod_{j=1}^n Y_{l_{j}\, j}
\ef,
\ee
because only $l_i$'s which are permutations in $\scrs_n$ have non-vanishing contribution in the sum.
This remark would be already enough to set the Cauchy--Binet theorem in
the simple case in which $X$ commutes with $Y$~\cite[Proposition~3.1]{uscapelli}.

The second step of the proof comes from analysing which terms do arise from
commuting the factor $Y_{l_1 1}$ to the position between 
$X_{\sigma(1) l_1}$ and $X_{\sigma(2) l_2}$, 
and so on recursively, by using the general formula
\be
x_1 [x_2\cdots x_r,y] = x_1 \sum_{s=2}^r x_2 \cdots x_{s-1} [x_s, y] x_{s+1} \cdots x_r
\ef.
\ee
As an illustration, we consider
the first application of this procedure
\be
\begin{split}
&
\bigg( \prod_{i=1}^n X_{\sigma(i) l_{i}} \bigg)
\bigg( \prod_{j=1}^n Y_{l_{j} j} \bigg)
=
X_{\sigma(1) l_{1}}
Y_{l_{1} 1}
\bigg( \prod_{i=2}^n X_{\sigma(i) l_{i}} \bigg)
\prod_{j=2}^nY_{l_{j} j}
\\
& \hspace{2cm}
+ \sum_{k=2}^n
\bigg( \prod_{r=1}^{k-1} X_{\sigma(r) l_{r}} \bigg)
[X_{\sigma(k) l_{k}},Y_{l_{1} 1}]
\bigg( \prod_{i=k+1}^n X_{\sigma(i) l_{i}} \bigg)
\prod_{j=2}^n Y_{l_{j} j}
\ef.
\end{split}
\ee
Then
\be
\label{eq.765876}
\begin{split}
& \sum_{\sigma \in \scrs_n} \sgn(\sigma) 
\sum_{l_{1},\dots, l_{n}\in [m]}
\bigg( \prod_{i=1}^n X_{\sigma(i) l_{i}} \bigg)
\prod_{j=1}^n Y_{l_{j} j}
\\
& =  \sum_{\sigma \in \scrs_n} \sgn(\sigma)  
\Bigg[
(XY)_{\sigma(1) 1}
\sum_{l_{2},\dots, l_{n}\in [m]}
\bigg( \prod_{i=2}^n X_{\sigma(i) l_{i}} \bigg)
\prod_{j=2}^n Y_{l_{j} j}
\\
& \hphantom{= \sum} 
- \sum_{k=2}^m
\sum_{l_{1},\dots, l_{n}\in [m]}
\bigg( \prod_{r=1}^{k-1} X_{\sigma(r) l_{r}} \bigg)
A_{\sigma(k) 1} B_{l_1 l_k}
\bigg( \prod_{i=k+1}^n X_{\sigma(i) l_{i}} \bigg)
\prod_{j=2}^n Y_{l_{j} j}
\Bigg]
\ef.
\end{split}
\ee
Consider the summands for each $k$ in the second row on the right hand side of
(\ref{eq.765876}). First of all, consider Lemma~\ref{lem.XXcomm} applied to a
matrix $X'$, defined as $X'_{ij} = X_{ij}$ if $i \neq k$ and
$A_{ij}$ if $i=k$. We are in the hypothesis of the Lemma because $X$
is row-pseudo-commutative and satisfies the condition~(\ref{XA2}). 
One can then write those summands as
\be
-
\sum_{\sigma \in \scrs_n} \sgn(\sigma) 
\!\!\!
\sum_{l_{1},\dots, l_{n}\in [m]}
\!\!\!
A_{\sigma(k) 1}
\bigg( \prod_{r=2}^{k-1} X_{\sigma(r) l_{r}} \bigg)
X_{\sigma(1) l_{1}}
B_{l_1 l_k}
\bigg( \prod_{i=k+1}^n X_{\sigma(i) l_{i}} \bigg)
\prod_{j=2}^n Y_{l_{j} j}
\ef.
\ee
Then, using the Translation Lemma for $\sigma \to \sigma \circ (1\,k)$,
and performing the sum over $l_1$
\be
+
\sum_{\sigma \in \scrs_n} \sgn(\sigma) 
\!\!\!
\sum_{l_{2},\dots, l_{n}\in [m]}
\!\!\!
A_{\sigma(1) 1}
\bigg( \prod_{r=2}^{k-1} X_{\sigma(r) l_{r}} \bigg)
(XB)_{\sigma(k) l_k} \\
\bigg( \prod_{i=k+1}^n X_{\sigma(i) l_{i}} \bigg)
\prod_{j=2}^n Y_{l_{j} j} 
\ef.
\ee
When $B_{ij} = \delta_{ij}$ the product of matrices $X$ becomes of the
same form of the first term of the right hand side of
(\ref{eq.765876}). This procedure can be repeated iteratively and,
ultimately, was enough to prove Proposition~\ref{theo.capCB}.

However, as the commutation of $X$'s and $Y$'s now
produces extra matrices $B$, we have to deal with an induction
expression of a more general form. One easily sees that, at all steps,
matrices $B$ will only act on $X$'s from the right, so, in order to
deal with the generic step $t$ of the procedure (beside $t=1$ seen in
detail above), we will consider expressions of the form
\be
\sum_{\sigma \in \scrs_n} \sgn(\sigma) L(\sigma)
\sum_{l_t\in [m]}
\bigg(\prod_{i=t}^n (X B^{\nu(i)})_{\sigma(i) l_{i}}\bigg)
\prod_{j=t}^nY_{l_{j} j}
\ee
where $L(\sigma)$ depend only from $\sigma_1,\dots\sigma_{t-1}$ and
$\nu(i)$ are non-negative integers. This form includes the initial
situation at $t=0$, and, as we see in a moment, is stable when $t$ is
increased. Indeed we have
\be
\begin{split}
\label{eq.765876d}
&
\sum_{\sigma \in \scrs_n} \sgn(\sigma) L(\sigma)
\sum_{l_t \in [m]}
\bigg(\prod_{i=t}^n (X B^{\nu(i)})_{\sigma(i) l_{i}}\bigg)
\prod_{j=t}^n Y_{l_{j} j}
\\
& \qquad
=
\sum_{\sigma \in \scrs_n} \sgn(\sigma) L(\sigma)
(X B^{\nu(t)} Y)_{\sigma(t) t}
\bigg(\prod_{i=t+1}^n (X B^{\nu(i)})_{\sigma(i) l_{i}}\bigg)
\prod_{j=t+1}^n Y_{l_{j} j}
\\
& \qquad \quad
+
\sum_{k=t+1}^n
\sum_{\sigma \in \scrs_n} \sgn(\sigma) L(\sigma)
\sum_{l_t\in [m]}
\bigg(\prod_{r=t}^{k-1} (X B^{\nu(r)})_{\sigma(r) l_{r}}\bigg)
\\
& \qquad \quad \qquad
\times \ 
A_{\sigma(k) t} (B^{\nu(k)+1})_{l_t l_k}
\bigg(\prod_{i=k+1}^n (X B^{\nu(i)})_{\sigma(i) l_{i}}\bigg)
\prod_{j=t+1}^n Y_{l_{j} j}
\ef.
\end{split}
\ee
In the last summands, we would like to commute the term
$A_{\sigma(k)t}$ in front of all $X$'s, as it carries the smallest
column-index.  This is indeed possible, at the light of
Lemma~\ref{lem.XXcomm}. Consider this lemma applied to a matrix $X'$,
defined as $X'_{ij} = (XB^{\nu(j)})_{ij}$ if $i \neq k$ and $A_{ij}$
if $i=k$. We are in the hypothesis of the Lemma because $X$ is
row-pseudo-commutative and satisfies the condition~(\ref{XA2}), and
therefore the same is true when replacing $X$ with $XB^{\nu(j)}$
because $B^{\nu(j)}$ acts on the column indices.  Then apply the
Involution Lemma with $(t\,k)$, and sum over $l_t$ where
appropriate. We can thus write
\be
\label{eq.765876f}
\begin{split}
&
\sum_{\sigma \in \scrs_n} \sgn(\sigma) L(\sigma)
\sum_{l_t\in [m]}
\bigg(\prod_{i=t}^n (X B^{\nu(i)})_{\sigma(i) l_{i}}\bigg)
\prod_{j=t}^n Y_{l_{j} j}
\\
& \qquad
=
\sum_{\sigma \in \scrs_n} \sgn(\sigma)  L(\sigma)\,
(X B^{\nu(t)} Y)_{\sigma(t) t}
\bigg(\prod_{i=t+1}^n (X B^{\nu(i)})_{\sigma(i) l_{i}}\bigg)
\prod_{j=t+1}^n Y_{l_{j} j} 
\\
&
\qquad \quad +
\sum_{k=t+1}^n
\sum_{\sigma \in \scrs_n} \sgn(\sigma)  L(\sigma)\,
A_{\sigma(t) t} 
\bigg(\prod_{r=t+1}^{k-1} (X B^{\nu(r)})_{\sigma(r) l_{r}}\bigg)
\\
&
\qquad \quad \qquad
\times \ 
(X B^{\nu(k)+\nu(t)+1})_{\sigma(k) k}
\bigg(\prod_{i=k+1}^n (X B^{\nu(i)})_{\sigma(i) l_{i}}\bigg)
\prod_{j=t+1}^n Y_{l_{j} j} 
\ef.
\end{split}
\ee
The relevant point in this expression is that all of the $n-t+1$
summands are of the same form of the original left hand side, with one
less matrix $Y$ to be reordered. However, while in the simpler case
$B_{ij}=\delta_{ij}$ the various terms were \emph{identical} up to the
prefactor, and could be collected together in a simple induction, here
they differ in the set of exponents $\{ \nu(i) \}$. Not accidentally,
the combinatorics of these lists of exponents has already been
discussed in Section~\ref{sec.enumeration}.
Indeed we can identify 
\begin{multline}
\label{eq.symb}
\nunu{\nu_1, \cdots, \nu_{t-1} | \nu_t, \cdots, \nu_n }
:= \\
\sum_{\sigma \in \scrs_n} \sgn(\sigma) 
\bigg( \prod_{i=1}^{t-1} M^{(\nu_{i})}_{\sigma(i) i} \bigg)
\sum_{l_t, \ldots, l_n\in [m]}
\bigg( \prod_{j=t}^n (X B^{\nu_{j}})_{\sigma(j) l_{j}} \bigg)
 \prod_{r=t}^n Y_{l_{r} r}
\ef,
\end{multline}
where parameters $\nu_i$ have to be integers, and $\nu_i \geq -1$ for 
$i=1, \ldots, t-1$, while $\nu_i \geq 0$ for $i=t, \ldots, n$. 
The matrix elements $M_{ij}^{(\nu_j)}$ are $A_{ij}$ if
$\nu_j=-1$ 
and $(X B^{\nu_j} Y)_{ij}$ if $\nu_j$ is non-negative.
In particular
\be
\vec{\nu}_{\bm 0} = \nunu{\, | \underbrace{0 \cdots 0}_n } = 
\sum_{\sigma \in \scrs_n} \sgn(\sigma)
\!
\sum_{l_{1},\dots, l_{n}\in [m]}
\!
\bigg( \prod_{i=1}^n X_{\sigma(i)\, l_{i}} \bigg)
 \prod_{j=1}^n Y_{l_{j}\, j}
\ee
Our rule (\ref{eq.765876f}) coincides with \reff{eq.symbrule} under
this identification, and we can apply Lemma~\ref{main} to get
\be
\label{open}
\sum_{\substack{ L \subseteq [m] \\ |L|=n }} 
\coldet X_{[n],L} \coldet Y_{L,[n]} 
=
\sum_\gamma c(\gamma) \sum_{\sigma \in \scrs_n} \sgn(\sigma) 
\bigg( \prod_{i=1}^{n} M^{(\nu_{i}(\gamma))}_{\sigma(i) i} \bigg)
\ef,
\ee
where notations are as in Section \ref{sec.enumeration}, i.e.\ 
$\gamma$ is a directed path in the upper half-plane starting from
the origin, the heights $(h_0, \dots, h_{t-1})$, are given by
$h_{i+1}-h_i=-\nu_i$, each $\nu_i$ is in the set $\{-1, 0, 1,
2,\ldots \}$, and the coefficients $c(\gamma)$ are given
by~\reff{eq.cgam}.

Now we can use Lemma~\ref{lemma:WH} to obtain
\be
\label{eq.154327652}
\begin{split}
&
\sum_{\substack{ L \subseteq [m] \\ |L|=n }} 
\coldet X_{[n],L} \coldet Y_{L,[n]}
= 
\sum_{\vec{\nu}_n} \bra{0 } \sum_{\sigma \in \scrs_n} \sgn(\sigma) 
\bigg( \prod_{i=1}^{n} \chi(\nu_i) \, M^{(\nu_{i})}_{\sigma(i) i} \bigg) \ket{0} 
\\
& 
\hspace{3cm}
= \bra{0}
\sum_{\sigma \in \scrs_n} \sgn(\sigma) 
\prod_{i=1}^{n} \bigg( \sum_{\nu_i = -1}^\infty   \chi(\nu_i) \,
M^{(\nu_{i})} \bigg)_{\sigma(i) i} \ket{0} 
\\
&
\hspace{3cm}
= \bra{0} 
\coldet \bigg( \sum_{\nu = -1}^\infty \chi(\nu) \, M^{(\nu)} \bigg) 
\ket{0}
\ef,
\end{split}
\ee
but
\be
\begin{split}
\sum_{\nu = -1}^\infty   \chi(\nu) \, M^{(\nu)}  
&= 
a\, M^{(-1)} + \sum_{\nu = 0}^\infty \big(\adg\big)^\nu M^{(\nu)}
= a\, A + X  \sum_{\nu = 0}^\infty \big(\adg B\big)^\nu Y
\\
&= a\, A + X \big(I - \adg B\big)^{-1} Y
\ef,
\end{split}
\ee
so we got our thesis. \qed

\section{Holomorphic representation}
\label{sec.holom}

The results of Proposition~\ref{new.capCB} can also be expressed as a
multiple integral in the complex plain, a structure that, within the
language of the quantum oscillator, is called a {\em holomorphic
  representation}. We shall use the {\em coherent states} of the
quantum oscillator, which are the states $\ket{z}$ defined as
\be
\ket{z} := \exp ( z \adg) \ket{0}
\ee
with $z\in \mathbb{C}$ a complex number. From the commutation
relations~\reff{WH} it soon follows the fundamental property of these
states
\be
\aaa \ket{z} = z \ket{z}
\ee
that is, it is an eigenstate of the annihilation operator. And, of course
\begin{align}
\bra{z} 
&:=
\bra{0}\, \exp ( \bar{z} \aaa)
\ef,
&
\bra{z} \adg 
&=
\bra{z} \bar{z}
\ef,
\end{align}
where $\bar{z}$ is the complex-conjugate of $z$. One easily verifies
that two different coherent states are not orthogonal
\be
\braket{z}{z'}
= \exp (\bar{z} z')
\ef.
\ee
However, since coherent states obey a closure relation, any state can
be decomposed on the set of coherent states. They hence form an
overcomplete basis. This closure relation can be expressed by the
resolution of the identity
\be
\int \frac{dz\,d\bar{z}}{i \pi} \exp (- |z|^2) \ket{z}\bra{z} = 1
\ef.
\ee
Let us consider the evaluation of
\be
\bra{0} 
\big( f_1(\adg) + g_1(\aaa) \big)
\cdots
\big( f_n(\adg) + g_n(\aaa) \big)
\ket{0}
\ef,
\ee
where $\{ f_{\alpha}, g_{\alpha} \}_{1 \leq \alpha \leq n}$ are $2n$ generic
expressions in a ring $R$, for which we have \emph{a priori} no
knowledge on the commutators\footnote{We mean here that, for
  $f(\adg)=\sum_i (\adg)^i f_i$, $g(\aaa)=\sum_j \aaa^j g_j$, with
  $f_i$'s and $g_j$'s in a \emph{commutative} ring,
  $[f(\adg),g(\aaa)]=\sum_{i,j} f_i g_j [(\adg)^i,\aaa^j]$, and the
  commutators are known, although complicated in general. However, if
  the coefficients $f_i$'s and $g_j$'s are valued in a generic
  \emph{non-commutative} ring, even if commuting with the
  Weyl--Heisenberg algebra, we have unknown extra terms of type
  $[f_i,g_j]$, namely: $[f(\adg),g(\aaa)]=\sum_{i,j} \big( g_j f_i
  [(\adg)^i,\aaa^j] + (\adg)^i \aaa^j [f_i, g_j] \big)$.}.  We are
ultimately interested in the case, corresponding to
Proposition~\ref{new.capCB},
\be
\bra{0} \coldet \big( F(\adg) + G(\aaa) \big) \ket{0} 
=
\sum_{\sigma  \in \scrs_n} 
\sgn(\sigma) 
\bra{0} 
\prod_{j=1}^n 
\big(F_{\sigma(j) j}(\adg) + G_{\sigma(j) j}(\aaa)\big) \ket{0} 
\ef,
\ee
(the product is ordered), with
\begin{align}
F(\adg) 
&= X (1 - \adg B)^{-1} Y
\ef,
&
G(\aaa) 
&=
\aaa A 
\ef.
\end{align}
Let $z_0=z_n=0$, and introduce $n-1$ intermediate coherent states,
with parameters $z_1, \ldots, z_{n-1}$, to get (with no more need of
ordered products on the right hand side)
\be
\bra{0} 
\prod_{j=1}^n 
\big( f_j(\adg) + g_j(\aaa) \big)
\ket{0}
=
\int 
\prod_{j=1}^{n-1} 
\left(
\frac{dz_j\,d\bar{z}_j}{i \pi} 
e^{-\left|z_j\right|^2}
\right)
\prod_{j=1}^{n}
\bra{z_{j-1}} 
f_j(\adg) + g_j(\aaa)
\ket{z_j}
\ef.
\ee
Each scalar product is easily evaluated according to
\be
\bra{u} f(\adg) + g(\aaa) \ket{v}
=
\big( f(\bar{u}) + g(v) \big)\, 
e^{\bar{u} v}
\ef,
\ee
so that
\begin{align}
\begin{split}
&
\bra{0} 
\prod_{j=1}^n 
\big( f_j(\adg) + g_j(\aaa) \big)
\ket{0}
=
\int 
\prod_{j=1}^{n-1} 
\frac{dz_j\, d\bar{z}_j}{i \pi } 
\,
e^{-\sum_{j=1}^{n} \bar{z}_j (z_{j}-z_{j+1})}
\prod_{j=1}^n 
\big( f_j(\bar{z}_{j-1}) + g_j(z_j) \big)
\ef,
\end{split}
\end{align}
and in particular
\begin{align}
\begin{split}
\bra{0} \coldet \big(  \aaa A + X (1 - \adg B)^{-1} Y \big) \ket{0}
&=
\int 
\prod_{j=1}^{n-1} 
\frac{dz_j\, d\bar{z}_j}{i \pi } 
\,
e^{-\sum_{j=1}^{n} \bar{z}_j (z_{j}-z_{j+1})}
\coldet M(z)
\end{split}
\\
M_{ij}(z)
&=
A_{ij} \, z_j 
+
\big( X (1 - \bar{z}_{j-1} B)^{-1} Y \big)_{ij}
\ef.
\end{align}
The equation above, jointly with Proposition \ref{new.capCB}, provides
a representation of the non-commutative Cauchy--Binet expression in
terms of an integral over $n$ (commuting) complex variables. This
result is somewhat implicit in Proposition \ref{new.capCB}, and the
standard general facts on the holomorphic representation of the
quantum oscillator.

Let us however observe that, in Section \ref{sec.full}, we could have
derived \emph{directly} the holomorpic representation, from the
Cauchy--Binet left hand side, instead of the representation in terms
of creation and annihilation operators. We only need to follow a
different track at the very final step of the proof, where, in
equation (\ref{eq.154327652}), we use the combinatorial Lemma~\ref{lemma:WH}.

The equivalent lemma for coherent states is based on the
formula\footnote{Which is easily proven, e.g.\ in generating function,
\[
\sum_{p,q} \frac{\zeta^p \xi^q}{p! q!}
\int \frac{dz\,d\bar{z}}{i \pi}\, z^p\, \bar{z}^q
\exp \big( -\bar{z}(z - \eta) \big) = 
\int \frac{dz\,d\bar{z}}{i \pi}\, 
\exp \big( -\bar{z}(z - \eta) + \bar{z} \xi + \zeta z \big) = 
\exp \big( \zeta (\eta + \xi) \big)
\ef,
\]
while
\[
\sum_{p,q} \frac{\zeta^p \xi^q}{p! q!}
\frac{p!}{(p-q)!}\,\eta^{p-q}
=
\sum_{p,q} \frac{\zeta^p}{p!}
(\eta + \xi)^{p}
=
\exp \big( \zeta (\eta + \xi) \big)
\ef.
\]}
\be
\label{eq.intecapo}
\int \frac{dz\,d\bar{z}}{i \pi}\, z^p\, \bar{z}^q
\exp \big( - \bar{z}(z - \eta) \big) = 
\frac{p!}{(p-q)!}\,\eta^{p-q}
\ef,
\ee
and reads (using notations as described in
Section~\ref{sec.enumeration} for paths $\gamma$, symbols
$\vec{\nu}_n$, coefficients $c(\vec{\nu}_{n})$, and conversion between
$\nu_i$'s and $h_i$'s)
\begin{lemma}
\label{lemma:WHholo}
For $\nu \geq -1$, define the monomials
\be
\chi_i(\nu)= \begin{cases}
\bar{z}_{i-1}^{\nu} & 
\nu \geq 0 \ef;
\\
z_i & 
\nu = -1 \ef.
\end{cases}
\ee
Then, when the symbol $\vec{\nu}_n=\nunu{\nu_1, \cdots \nu_n | \,}$
corresponds to a path $\gamma$, setting $z_0=z_n=0$,
\be
\int 
\prod_{j=1}^{n-1}
\frac{dz_j\, d\bar{z}_j}{i \pi}
\,
e^{-\sum_{j=1}^{n-1} \bar{z}_j (z_{j}-z_{j+1})}
\chi_1(\nu_1) \cdots \chi_n(\nu_n)
= 
c(\vec{\nu}_{n})
=
\prod_{\substack{i \in [n] \\ h_i \leq h_{i-1} }} 
\frac{h_{i-1}!}{h_i!}
\ef,
\ee
while otherwise the integral
above is zero.
\end{lemma}

\proof We try to follow as closely as possible the reasoning in the
proof of Lemma \ref{lemma:WH}. We proceed by induction.  Assume that, for a sequence
$\nu_1,\dots, \nu_t$ such that the corresponding path remains
positive,
\be
\int 
\prod_{j=1}^{t-1}
\frac{dz_j\, d\bar{z}_j}{i \pi}
\,
e^{-\sum_{j=1}^{t-1} \bar{z}_j (z_{j}-z_{j+1})}
\chi_1(\nu_1) \cdots \chi_t(\nu_t) = 
z_t^{h_t}
\! \prod_{\substack{i \in [t] \\ h_i \leq h_{i-1} }} 
\!\! \frac{h_{i-1}!}{h_i!}
\ef.
\ee
This is indeed the case for $t=0$ (where, as customary for products
over empty sets, we have $1=1$), and in the more convincing case $t=1$
(where we have no integrations to perform, and, as $z_0=0$,
$\chi_1(\nu_1)= z_1$, $1$ and $0$ respectively if $\nu_1=-1$, $0$ or
strictly positive).


Then, we analyse the consequence of increasing $t$ on both sides of
the equation. On the left hand side, we should multiply by
$e^{-\bar{z}_t (z_{t}-z_{t+1})} \chi_{t+1}(\nu_{t+1})$, and then
integrate over $dz_t\, d\bar{z}_t$.  If $\nu_{t+1} = -1$,
$\chi_{t+1}(\nu_{t+1})=z_{t+1}$ and $h_{t+1} = h_t+1$,
while if $\nu_{t+1} \geq 0$,
$\chi_{t+1}(\nu_{t+1})=\bar{z}_{t}^{\nu_{t+1}}$ and
$h_{t+1} = h_t - \nu_{t+1}$.
In both cases, the integral is of the form (\ref{eq.intecapo}),
and we get
\begin{align}
\int 
\frac{dz_t\, d\bar{z}_t}{i \pi}
\,
e^{-\bar{z}_t (z_{t}-z_{t+1})}
\,
z_t^{h_t}
z_{t+1}
&=
z_{t+1}^{h_t + 1}
=
z_{t+1}^{h_{t+1}}
\ef;
\\
\int 
\frac{dz_t\, d\bar{z}_t}{i \pi}
\,
e^{-\bar{z}_t (z_{t}-z_{t+1})}
\,
z_t^{h_t}
\bar{z}_t^{\nu_{t+1}}
&=
\frac{h_t !}{(h_t - \nu_{t+1})!}
z_{t+1}^{h_t - \nu_{t+1}}
=
\frac{h_t !}{h_{t+1}!}
z_{t+1}^{h_{t+1}}
\ef.
\end{align}
In the two cases, the integration produces the
appropriate relative factor, which, in particular, is zero if the path
goes below the horizontal axis (because of a $1/k!$ factor, with
$k<0$). At the last step, we remain with a factor $z_n^{h_n}$. As
$z_n=0$, we select only the paths terminating at height zero. \qed

\section{A lemma on the Campbell-Baker-Hausdorff formula}
\label{sec.cbh}

The goal of this section is to prove the following relation, which is
a preparatory lemma to our Capelli identity in Grassmann
representation, proven in the next section.
\begin{prop}
\label{prop.cbh}
Let $\aaa$ and $\adg$ be the generators 
of a Weyl-Heisenberg Algebra, i.e.\ $[\aaa, \adg] = 1$,
and $f(x)$ a formal power series. Then,
at the level of formal power series, we have
\be
\label{eq.cbhinvres}
\begin{split}
\exp \left( \adg + f(\aaa) \right)
&= 
\exp(\adg) 
\exp \bigg( 
\sum_{k\geq 0}
\frac{1}{(k+1)!}
(\partial^k f)(\aaa) 
\bigg)
\\
&=
\exp(\adg)
\exp \bigg( 
\frac{\exp(\partial)-1}{\partial} f(\aaa)
\bigg)
\ef.
\end{split}
\ee
\end{prop}
The proposition above is a special case of the
\emph{Campbell-Baker-Hausdorff} (CBH) formula 
\cite{cbh_camp, cbh_poinc, cbh_baker, cbh_haus}.  We give here a proof
that makes use only of the existence of a CBH formula (and not the
explicit expressions known in the literature). Furthermore, an
additional argument provides a slightly longer variant, which instead
is completely self-contained.

We recall that, given two elements $x$ and $y$ in a non-commutative ring, the
Campbell-Baker-Hausdorff formula is an expression for
$\ln(\exp(x) \exp(y))$ as a formal infinite sum of elements of the
Lie algebra generated by $x$ and $y$:
\begin{align}
\exp(x) \exp(y) 
&= 
\exp \left(
x + y + z
\right)
\ef;
&
z
&=
S(x,y)
\ef;
\end{align}
The first few terms read
\be
S(y;x)=
\frac{1}{2} [x,y] + 
\frac{1}{12} [x-y,[x,y]] + \cdots
\ef,
\ee
and the generic summand in this series
has the form 
\[
[z_{s(1)},[z_{s(2)},\cdots [z_{s(k-1)},z_{s(k)}]\cdots ]]
\]
for some integer $k\geq 2$, $(s(1),\ldots,s(k)) \in \{0,1\}^k$, and the
identification $z_0=x$, $z_1=y$. Of course, terms with $s(k)=s(k-1)$
vanish in any Lie algebra, and many other strings are redundant, 
e.g., besides the trivial 
$[\cdots,[x,y]\cdots] = - [\cdots,[y,x]\cdots]$, 
a first non-trivial relation is
$[x,[y,[x,y]]]=[y,[x,[x,y]]]$.

The existence statement is relatively easy to obtain.  The full
expression at all orders with coefficients in closed form is
complicated, but redundant forms (in the sense above) are well-known
in the literature (see e.g.~\cite[pp.~134 and 135]{saglewalde}).

Formal inversion (that is, solving w.r.t.~$y$, leaving $z$ as an
indeterminate) is easily achieved.  Define the inverse problem as
\begin{align}
\label{eq.cbhinv}
\exp \left( x + z \right)
&= 
\exp(x) \exp(z + y) 
\ef;
&
y
&=
\tilde{S}(x;z)
\ef;
\end{align}
then, multiplying both sides by $e^{-x}$ from the left, one obtains
\be
\tilde{S}(x;z) = S(-x,x+z)
\ef.
\ee
The existence result for $\tilde{S}$ follows from existence for $S$
and the relation above.

\proofof{Proposition~\ref{prop.cbh}}
Our proposition corresponds to the solution of the inverse
problem~(\ref{eq.cbhinv}), finding an expression for $\tilde{S}(x;z)$,
in the special case of $x=\adg$ and $z=f(\aaa)$.

In this case many commutators vanish. We have
\be
[\underbrace{\adg,[\adg,\cdots [\adg}_{k}, f(\aaa)]\cdots]]
=
(-\partial)^k f(\aaa)
\ee
where $\partial^k f$ denotes the $k$-th derivative of $f$ (as a power
series). So, all the expressions above do commute with $f(\aaa)$ and
we see that in our case all non-vanishing strings are the ones of the
form $(0,0,\cdots,0,1)$ (the ones $(0,0,\cdots,0,1,0)$ are also
non-vanishing but clearly redundant). In other terms, writing for a
generic Lie algebra
\be
\tilde{S}(x;z)
=
\sum_{k \geq 1} c_k \, [\underbrace{x,[x,\cdots [x}_{k}, z]\cdots]]
\; + \mathcal{O}(z^2)
\ef,
\ee
(where $\mathcal{O}(\cdot)$ is in the sense of polynomials in the
enveloping algebra), we get in our case
\be
\label{eq.stilck}
\tilde{S}(\adg; f(\aaa))
=
\sum_{k \geq 1} c_k \, 
[\underbrace{\adg,[\adg,\cdots [\adg}_{k}, f(\aaa)]\cdots]]
=
\sum_{k} c_k \, (-\partial)^k
f(\aaa)
\ef.
\ee
Observe that, again in the enveloping algebra,
\be
\label{eq.65465765}
[\underbrace{x,[x,\cdots [x}_{k}, z]\cdots]]
=
\sum_{h=0}^{k} (-1)^h \binom{k}{h}
x^{k-h} z x^h
\ee
and that
\be
\exp(x+z)
=
\exp(x) +
\sum_{k \geq 0}
\sum_{h=0}^k
\frac{1}{(k+1)!}
x^{k-h} z x^h
+ \mathcal{O}(z^2)
\ef.
\ee
Appealing to the existence of a solution, we can determine the $c_k$'s
by matching the coefficient of $z x^k$ on the two sides of
(\ref{eq.cbhinv}), using (\ref{eq.stilck}) and (\ref{eq.65465765}),
obtaining
\be
c_k = \frac{1}{(k+1)!}
\ee
that, with the fact 
$\sum_{k \geq 0} x^k/(k+1)! = (e^{x}-1)/x$ (used here at the
level of formal power series), gives our statement.


Avoiding to appeal to the existence statement requires to match 
all possible other linear monomials, of the kind $x^h z x^{k-h}$.
Then, the consistency of the assignment of $c_k$'s boils down to the
following relation: for each $k$ and $h$ positive integers,
\be
\sum_{i=0}^h
(-1)^{h-i} 
\binom{k+1}{i}
\binom{k-i}{h-i}
= 1
\ef.
\ee
This is proven by observing that 
$\binom{k-i}{h-i} = (-1)^{h-i} \binom{h-k-1}{h-i}$, and using
Chu-Vandermonde convolution, 
$\sum_i \binom{n}{i} \binom{m}{k-i} = \binom{n+m}{k}$.
\qed

If instead of $\adg$ we have $c\, \adg$, with $c$ some commuting
quantity, the same reasoning can be done, and a simple scaling applies
to all formulas. The corresponding generalization of
(\ref{eq.cbhinvres}) is
\be
\label{eq.cbhinvresC}
\begin{split}
\exp \left( c\, \adg + f(\aaa) \right)
&= 
\exp(c\, \adg) 
\exp \bigg( 
\sum_{k\geq 0}
\frac{c^k}{(k+1)!}
(\partial^k f)(\aaa) 
\bigg)
\\
&=
\exp(c\, \adg)
\exp \bigg( 
\frac{\exp(c\, \partial)-1}{c\, \partial} f(\aaa)
\bigg)
\ef.
\end{split}
\ee
We shall need also the identity obtained by Hermitian conjugation
\be
\label{eq.cbhinvresCH}
\begin{split}
\exp \left( c\, \aaa + f(\adg) \right)
&= 
\exp \bigg( 
\sum_{k\geq 0}
\frac{c^k}{(k+1)!}
(\partial^k f)(\adg) 
\bigg) \exp(c\, \aaa) 
\\
&=
\exp \bigg( 
\frac{\exp(c\, \partial)-1}{c\, \partial} f(\adg)
\bigg) \exp(c\, \aaa)
\ef.
\end{split}
\ee

\section{The Capelli identity in Grassmann Algebra}
\label{sec.grass}

Besides column- and row-determinants, defined in \reff{def.coldet} and
\reff{def.rowdet} respectively,
another possible non-commutative generalization of the determinant
is the \emph{symmetric-determinant}: 
\be
\label{def.symdet}
\symdet M 
\bydef
\frac{1}{n!} 
\sum_{\sigma, \tau \in \scrs} \sgn(\sigma) \sgn(\tau)  \, \prod_{i=1}^n M_{\sigma(i) \tau(i)} 
\ef.
\ee
In contrast to the cases of the column- and row-determinant, the
definition~\reff{def.symdet} demands in general the inclusion of
rational numbers in the field $K$ over which the ring $R$ is defined.

For any permutation $\tau \in \scrs_n$ let us denote $M^{\tau}$ the
matrix with entries $(M^{\tau})_{ij} = M_{i \, \tau(j)}$, and
${}^{\tau}M$ the matrix with entries 
$({}^{\tau}M)_{ij} = M_{\tau(i)\, j}$.  We clearly have, for any
matrix $M$, 
\begin{align}
\coldet {}^{\tau}M 
&=
\sgn(\tau)  \coldet M
\ef;
&
\rowdet M^{\tau} 
&=
\sgn(\tau) \rowdet M
\ef;
\end{align}
while in general the action of the symmetric group on columns and
rows, respectively for the two cases, is not simple.

Indeed, the symmetric-determinant reads
\be
\label{id.symdet}
\symdet M 
=
\frac{1}{n!}
\sum_{\tau  \in \scrs_n} \sgn(\tau) 
\coldet M^{\tau}
=
\frac{1}{n!}
\sum_{\tau  \in \scrs_n} \sgn(\tau) 
\rowdet {}^{\tau}M
\ef,
\ee
and no relevant further simplifications are possible in general.

However, for a $n$-dimensional matrix $M$ with weakly row-symmetric
commutators, (and thus in particular if $M$ is row-pseudo-commutative),
in~\cite[Lemma 2.6(a)]{uscapelli} we proved that
\emph{both} actions of the symmetric group are simple, i.e.\ also
\begin{align}
\coldet M^{\tau} 
&=
\sgn(\tau)  \coldet M 
\ef;
\label{id.coldet}
\end{align}
(and similarly for the row-determinant, if $M$ has
weakly column-symmetric commutators),
and therefore for such a matrix the expression (\ref{id.symdet}) simplifies
(in particular, rationals are not necessary)
\begin{coroll}
For a $n$-dimensional matrix $M$ with weakly row-symmetric  commutators
\be
 \symdet M = \coldet M 
\ef.
\ee
\end{coroll}
Our interest in the symmetric-determinant follows from the remark that
it provides the generalization of the Berezin integral
representation~\reff{simplegrass} for the determinant of a matrix with
commuting elements.  Indeed, for $M$ a $n \times n$ matrix with
elements in a non-commutative ring $R$, if $R$ contain the rationals
(or $M$ is row-pseudo-commutative), and $\{\psibar_i, \psi_i\}_{i
  \in [n]}$ a set of $2n$ Grassmann variables commuting with the
entries $M_{ij}$, we have
\be
\label{symdetint}
\int \mathcal{D}(\psi, \psibar)
\exp(\psibar M \psi)
=
\symdet M
\ef.
\ee
Comparatively, the Grassmann formulas for the
column- and row-determinant are more cumbersome, as they require an
ordering of the $n$ factors
\begin{align}
\label{coldetint}
\int d \psi_n \cdots d \psi_1
\;
(\psi M)_1 \cdots (\psi M)_n
&=
\coldet M
\ef;
\\
\label{rowdetint}
\int d \psi_n \cdots d \psi_1
\;
(M \psi)_1 \cdots (M \psi)_n
&=
\rowdet M
\ef.
\end{align}
Grassmann indeterminates present the advantage of encoding
our commutation relations in a simple way.  For example:
\begin{lemma}
\label{lemma:C}
Let $R$ be a 
ring, and $A$ a $n\times n$ matrix with elements in $R$. Let the
$\{\psi_i\}_{i \in [n]}$ be nilpotent Grassmann indeterminates, that
is $\psi_i^2 =0$ and their anti-commutators $\{\psi_i, \psi_j\}=0$
vanish.
\begin{itemize}
\item[(a)] Let $X$  be a $n\times m$  matrix with elements in $R$ such that
\be
\label{eq.1547326543}
[X_{ij}, A_{k\ell}] - [X_{kj}, A_{i\ell}] = 0 \qquad\hbox{\rm for all } i,j,k,\ell 
\ef.
\ee
then
\be
\{ (\psi X)_j, (\psi A)_{\ell} \} := 
\sum_{i\in [n]} \sum_{k \in [n]} \{
\psi_i X_{ij}, \psi_k A_{k \ell} \} = 0
\ef.
\ee
\item[(b)] Let $Y$  be a $m\times n$  matrix with elements in $R$ such
  that
\be
[Y_{ij}, A_{k\ell}] - [Y_{i\ell}, A_{kj}] = 0 
\qquad\hbox{\rm for all } i,j,k,\ell 
\ef.
\ee
then
\be
\{ (Y \psi )_i, (A \psi )_k \} := \sum_{j\in [n]} \sum_{\ell \in [n]} 
\{ Y_{ij} \psi_j , A_{k \ell} \psi_{\ell} \} = 0
\ef.
\ee
\end{itemize}
\end{lemma}

\proof
(a) We have that
\be
\begin{split}
\{ (\psi X)_j, (\psi A)_{\ell} \} 
&=
\sum_{i, k \in [n]}
\left( \psi_i X_{ij} \psi_k A_{k\ell} + \psi_k A_{k\ell} \psi_i X_{ij}
\right)
=
\sum_{i, k \in [n]}
\psi_i \psi_k \, [ X_{ij}, A_{k\ell} ]
\\
& =
\sum_{1 \leq i < k \leq n}
\psi_i \psi_k 
\,\Big( [ X_{ij}, A_{k\ell} ] - [ X_{kj}, A_{i\ell} ] \Big)
\ef,
\end{split}
\ee 
where we have taken into account that $i\neq k$ because the $\psi$'s
are nilpotent and we have put together the terms in which both
$\psi_i$ and $\psi_k$ appears. But now each term in the sum vanish by
the hypothesis~(\ref{eq.1547326543}).
The case (b) is identical.
\qed

This result is used to prove the following:
\begin{lemma}
\label{commute}
Let $R$ be a
ring, and $X$ a $n\times m$, $Y$ a $m\times n$, $A$ a $n\times n$ and
$B$ a $m\times m$ matrix with elements in $R$. Let the $\{\psibar_i,
\psi_i\}_{i \in [n]}$ be nilpotent Grassmann indeterminates commuting
with $R$, that is $\psibar_i^2=\psi_i^2 =0$ and their anti-commutators
$\{\psibar_i, \psibar_j\}=\{\psibar_i, \psi_j\}=\{\psi_i, \psi_j\}=0$
vanish.  If
\be
[X_{ij}, A_{k\ell}] - [X_{kj}, A_{i\ell}] = [Y_{ij}, A_{k\ell}] -
[Y_{i\ell}, A_{kj}] = 0 \qquad\hbox{\rm for all } i,j,k,\ell 
\ee
and the elements of $B$ commute with the ones of $A$,
then for each integer $s$
\be
[ \psibar X B^s Y \psi, \psibar A \psi ]
= 0
\ef.
\ee
\end{lemma}

\proof Indeed, as $B_{ij}$'s and $A_{k \ell}$'s do commute, we can
write the commutator as
\be
[ \psibar X B^s Y \psi, \psibar A \psi ] =  
\sum_{r\in [n]} 
(\psibar X B^s)_{r} [ (Y \psi)_r, \psibar A \psi] 
+
[(\psibar X)_r , \psibar A \psi]\,(B^s Y \psi)_r 
\ef.
\ee
Consider separately each of the resulting commutators:
\begin{align}
\label{eq.64375658u}
[(\psibar X)_r , \psibar A \psi] 
&=
\sum_{k\in [n]} 
\{(\psibar X)_r, (\psibar A)_k\}\, \psi_k  = 0 
\ef;
\\
[ (Y \psi)_r, \psibar A \psi] 
&=
\sum_{k\in [n]}  
\psibar_k \, \{ (Y \psi)_r, (A \psi)_k\} =0
\ef;
\end{align}
where we used Lemma~\ref{lemma:C}.
\qed

We have now all the ingredients to prove Proposition~\ref{prop:grass}.

\proofof{Proposition~\ref{prop:grass}} (a) 
As $Y$ is row-pseudo-commutative, and
we assumed that our ring contains the rationals,
using (\ref{id.coldet}),
we can rewrite the left hand side of (\ref{eq.inPropGrass}) as
\be 
\label{eq.6476576}
\sum_{\substack{L \subseteq [m] \\ |L|= n}}
\coldet X_{[n],L}  \coldet Y_{L,[n]}
=
\frac{1}{n!}
\sum_{\tau  \in \scrs_n}
\sgn(\tau)
\sum_{\substack{L \subseteq [m] \\ |L|= n}}
\coldet X_{[n],L} \coldet (Y^{\tau})_{L,[n]} 
\ef.
\ee
From the hypotheses we soon have that, for any permutation $\tau \in
\scrs_n$, the matrices $X, Y^{\tau}, A^{\tau}, B$ satisfy the
hypothesis of Proposition~\ref{new.capCB}(a), and therefore, as $X$ is
row-pseudo-commutative, we have that
\be
\sum_{\substack{L \subseteq [m] \\ |L|= n}}
\coldet X_{[n],L} \coldet (Y^{\tau})_{L,[n]} =
\bra{0}
\coldet (\aaa A^{\tau} + X (1 - \adg B)^{-1} Y^{\tau})
\ket{0}
\ef.
\ee
Note that, on the right hand side, the permutation $\tau$ has exactly
the action from the right on the matrix $M = \aaa A + X (1 - \adg
B)^{-1} Y$. Thus, the combination in (\ref{eq.6476576}) corresponds to
the definition~\reff{id.symdet} of the symmetric-determinant,
\be
\frac{1}{n!}
\sum_{\tau \in \scrs_n } 
\sgn(\tau)
\coldet (\aaa A^{\tau} + X (1 - \adg B)^{-1} Y^{\tau})
=
\symdet (\aaa A + X (1 - \adg B)^{-1} Y)
\ef.
\ee
We can use the Grassmann representation, \reff{symdetint}, for the
expression above,
to conclude that
\be
\sum_{\substack{L \subseteq [m] \\ |L|= n}}
\coldet X_{[n],L}  \coldet Y_{L,[n]} = \int \mathcal{D}(\psi, \psibar)
\bra{0}
\exp \big(
\psibar A \psi\,\aaa  +  \psibar X (1 - \adg B)^{-1} Y \psi
\big)
\ket{0}
\ef.
\ee
Now we use the result in (\ref{eq.cbhinvresCH}) by posing
$c = \psibar A \psi$
and
$f(\adg) = \psibar X (1 - \adg B)^{-1} Y \psi$.
Using the hypotheses (\ref{eq.XAgrass}) and (\ref{eq.YAgrass}) of the
proposition, we can verify the hypothesis of Lemma~\ref{commute},
therefore our quantities $c$ and $f(\adg)$ commute (as required for
(\ref{eq.cbhinvresCH}) to apply), and we get
\be
\label{159}
\exp \big(
\psibar A \psi\,\aaa  +  \psibar X (1 - \adg B)^{-1} Y \psi \big) = \exp( g(\adg) ) \exp(
\psibar A \psi \, \aaa )
\ee
with $g(\adg)$ determined according to
(\ref{eq.cbhinvresCH}),\footnote{Note at this aim
  that, if $[M_{ij},M_{k\ell}]=0$, 
  $\frac{\partial}{\partial \xi} (\vec{u}, (I-\xi M)^{-s} \vec{v})
   = s \, (\vec{u} M, (I-\xi M)^{-s-1} \vec{v})$.}
\be
g(\adg)
= 
\sum_{k \geq 0}
\frac{(\psibar A \psi)^{k}}{k+1}
( \psibar X B^k (1 - \adg B)^{-k-1} Y \psi )
\ef.
\ee
Note that, in the sum, $k$ cannot become larger that $n-1$, because of the
nilpotency of the Grassmann indeterminates.

In~\reff{159} the creation and annihilation operators are
ordered into a polynomial with monomials of the form $(\adg)^k \aaa^h$
(i.e., they are \emph{antinormal--}, or {\em anti-Wick--ordered}),
and the whole expression is drastically simplified because
\begin{align}
\exp( \aaa \, \psibar A \psi ) \, \ket{0}
& =
\ket{0}
\ef;
\\
\bra{0} \exp( g(\adg) )
=
\bra{0}\exp( g(0) )
&=
\bra{0}
\exp 
\bigg(
\sum_{k \geq 0}
\frac{(\psibar A \psi)^{k}}{k+1}
( \psibar X B^k Y \psi )
\bigg)
\ef.
\end{align}
As there are no more creation and annihilation operators,
we can just drop the factor $\braket{0}{0}=1$, to obtain the purely
fermionic representation
\be
\sum_{\substack{L \subseteq [m] \\ |L|= n}}
\coldet X_{[n],L}  \coldet Y_{L,[n]}
=
\int \mathcal{D}(\psi, \psibar)
\exp 
\bigg(
\sum_{k \geq 0}
\frac{(\psibar A \psi)^{k}}{k+1}
(\psibar X B^k Y \psi)
\bigg)
\ef,
\ee
or, by summing over $k$, intending $\ln (I-M) = \sum_{k \geq 1}
\frac{1}{k} M^k$ as a polynomial, truncated by the nilpotence of
$\psibar A \psi$, and using 
$[ (\psibar X)_r, \psibar A \psi ] = 0$ for every $r$
(valid because of Lemma \ref{lemma:C}, see equation
(\ref{eq.64375658u})),
\be
\sum_{\substack{L \subseteq [m] \\ |L|= n}}
\coldet X_{[n],L}  \coldet Y_{L,[n]}
=
\int \mathcal{D}(\psi, \psibar)
\exp \bigg(
- \psibar X
\frac{\ln (1 - (\psibar A \psi) B)}{(\psibar A \psi) B}
Y \psi
\bigg)
\ef,
\ee
as announced.

For the case (b), consider now the matrices ${}^{\tau}X, Y, {}^{\tau}A, B$ which
satisfy the hypothesis of Proposition~\ref{new.capCB}(b) and
therefore, as $X$ and $Y$ are column-pseudo-commutative,
following the procedure above,
\be
\begin{split}
&
\sum_{\substack{L \subseteq [m] \\ |L|= n}}
\!
\rowdet X_{[n],L} \rowdet Y_{L,[n]}
= \int \mathcal{D}(\psi, \psibar)
\bra{0}
\exp \big(
\adg \psibar A \psi  +  \psibar X (1 - \aaa B)^{-1} Y \psi
\big)
\ket{0}
\end{split}
\ee
and, to conclude, we proceed as in the previous case, except that we use
the identity~(\ref{eq.cbhinvresC}) instead of 
(\ref{eq.cbhinvresCH}).
\qed


\section{Direct proof of the Grassmann representation for $B=I$}
\label{sec.directGrass}

We have proven a Grassmann version of the non-commutative
Cauchy--Binet formula as a consequence of the Weyl--Heisenberg
version. Considering also the necessary analysis of combinatorics of
{\L}ukasiewicz paths, for the latter, and of Campbell-Baker-Hausdorff
formula, for the former, the proof is quite composite.  It is
conceivable that a more direct proof may exist.

In this section we give such a proof, in the simplified situation in
which, besides the hypotheses in Proposition \ref{prop:grass}, we have
that $B$ is the identity matrix.  Indeed, in this case, the version of
non-commutative Cauchy--Binet formula obtained in \cite{uscapelli}
(and reported here as Proposition \ref{theo.capCB}(a)), and the
Grassmann-Algebra representation of Proposition \ref{prop:grass}(a),
hold simultaneously.  We produce here a short proof of the specialized
Proposition \ref{prop:grass}(a), taking Proposition
\ref{theo.capCB}(a) as the starting point.

Actually, just like in Proposition \ref{prop.old}, we will end up
proving that this relation between the right hand sides of
(\ref{eq.thCapCB}) and (\ref{eq.inPropGrass}) is in fact valid
regardless from the fact that $A$ is related to the commutator of $X$
and $Y$, i.e.\ they are a consequence
of a stronger fact
\begin{prop}
\label{prop.comeP13perGrass}
Let $R$ be a ring containing the rationals, and $U$ and $V$ be two
$n\times n$ matrices with elements in $R$. Let $\psibar_i$, $\psi_i$,
with $1 \leq i \leq n$, be Grassmann indeterminates. 
Define
\be
\label{eq.defQcol_rep9}
\big( Q^{\rm col}(V) \big)_{ij} \bydef V_{i j} (n-j)
\ef.
\ee
Assume that
\be
\label{eq.789765454}
[\psibar U \psi, \psibar V \psi]=0
\ef,
\ee
and that, for any permutation $\tau$,
\be
\label{eq.symcdet91}
\sgn(\tau)
\coldet (U^{\tau} + Q^{\rm col}(V^{\tau})) 
=
\coldet (U + Q^{\rm col}(V)) 
\ef.
\ee
Then
\be
\label{eq.inP91}
\coldet (U + Q^{\rm col}(V)) 
=
\int \mathcal{D}(\psi, \psibar)
\exp 
\bigg(
\sum_{k \geq 0}
\frac{(\psibar V \psi)^{k}}{k+1}
(\psibar U \psi) 
\bigg)
\ef.
\ee
\end{prop}

\proof
Remark that, for $s$ and $t$ commuting indeterminates, at the level of
power series,
\be
\exp 
\bigg(
s \sum_{k\ge 0} \frac{t^{k+1}}{k+1}
\bigg)
=
(1-t)^{-s}
= \sum_{n\ge 0}\,  \frac{t^n}{n!}\,  \big(s + (n-1)\big) \big(s + (n-2)\big) \dots s
\ef.
\ee
With the choice
$t\to t v$ and $s\to u/(t v)$, with
$u$, $v$ and $t$ commuting,
we get that
\be
\label{eq.453432654}
\exp
\bigg(
t\, u \sum_{k\ge 0} \frac{(t\,v)^k}{k+1} 
\bigg)
=
\sum_{n\ge 0}\, \frac{t^n}{n!}\, 
\big(u + (n-1) v\big) \big(u + (n-2) v\big) \dots u
\ef.
\ee
We apply this formula to the right hand side of (\ref{eq.inP91}), with
$u = \psibar U \psi$, $v = \psibar V \psi$, and and $t$ a formal
indeterminate that counts the degree in Grassmann variables (the
coefficient of order $t^k$ has $k$ factors $\psibar_i$'s and $k$
$\psi_j$'s). In particular, Grassmann integration selects only the
term $t^n$, and we get
\begin{multline}
\label{eq.43254464}
\int \mathcal{D}(\psi, \psibar)
\exp 
\bigg(
\sum_{k \geq 0}
\frac{(\psibar V \psi)^{k}}{k+1}
(\psibar U \psi) 
\bigg)
\\
=
\frac{1}{n!}\, \int \mathcal{D}(\psi, \psibar)
\,\big( \bar{\psi}\big( U + V (n-1)\big)\psi\big)
\,\big( \bar{\psi}\big( U + V (n-2)\big)\psi\big)
\, \dots \big(\bar{\psi}U\psi\big)
\ef.
\end{multline}
The left hand side of (\ref{eq.inP91}), using (\ref{coldetint}), reads
\be
\label{eq.41434354}
\int d \psibar_n \cdots d \psibar_1
\big( \psibar ( U +Q^{\rm col} ) \big)_1
\cdots
\big( \psibar ( U +Q^{\rm col} ) \big)_n
\ef,
\ee
that is, given the expression (\ref{eq.defQcol_rep9}) for $Q^{\rm col}$,
\be
\label{eq.6535654}
\int d \psibar_n \cdots d \psibar_1
\big( \psibar ( U + V (n-1) ) \big)_1
\big( \psibar ( U + V (n-2) ) \big)_2
\cdots
\big( \psibar U \big)_n
\ef.
\ee
We can introduce a trivial factor
$
1
=
\int d \psi_n \cdots d \psi_1
\psi_1 \cdots \psi_n
$,
and reorder the Grassmann variables, and terms in the integration
measure, to rewrite (\ref{eq.6535654}) as
\be
\label{eq.6535654b}
\int 
\mathcal{D}(\psi, \psibar)
\big( \psibar ( U + A (n-1) ) \big)_1
\big( \psibar ( U + A (n-2) ) \big)_2
\cdots
\big( \psibar U \big)_n
\psi_n \cdots \psi_1
\ef.
\ee
We can exploit the invariance in the hypothesis (\ref{eq.symcdet91}),
and the fact that our ring contains the rationals, to replace the
expression above by its symmetrization
\be
\label{eq.6535654btausym}
\frac{1}{n!}
\sum_{\tau}
\sgn(\tau)
\int 
\mathcal{D}(\psi, \psibar)
\big( \psibar ( U^{\tau} + A^{\tau} (n-1) ) \big)_1
\cdots
\big( \psibar U^{\tau} \big)_n
\psi_n \cdots \psi_1
\ef.
\ee
As $(M^{\tau})_{ij} = M_{i \tau(j)}$, we just have
\be
\label{eq.6535654btausym2}
\frac{1}{n!}
\sum_{\tau}
\sgn(\tau)
\int 
\mathcal{D}(\psi, \psibar)
\big( \psibar ( U + A (n-1) ) \big)_{\tau(1)}
\cdots
\big( \psibar U \big)_{\tau(n)}
\psi_n \cdots \psi_1
\ef.
\ee
Note that the factors $(n-j)$, multiplying the matrix entries of $A$,
remain unchanged in their ordering, and in particular the values of
$j$ are distinct from the indices, now $\tau(j)$, in the corresponding
product. Reorder the factors $\psi_i$'s so to compensate for the
signature of the permutation
\be
\label{eq.6535654btausym3}
\frac{1}{n!}
\sum_{\tau}
\int 
\mathcal{D}(\psi, \psibar)
\big( \psibar ( U + A (n-1) ) \big)_{\tau(1)}
\cdots
\big( \psibar U \big)_{\tau(n)}
\psi_{\tau(n)} \cdots \psi_{\tau(1)}
\ef,
\ee
and extend the sum to all $n$-uples of integers
\be
\label{eq.6535654btausym4}
\frac{1}{n!}
\sum_{i_1, \ldots, i_n \in [n]}
\int 
\mathcal{D}(\psi, \psibar)
\big( \psibar ( U + A (n-1) ) \big)_{i_1}
\big( \psibar ( U + A (n-2) ) \big)_{i_2}
\cdots
\big( \psibar U \big)_{i_n}
\psi_{i_n} \cdots \psi_{i_1}
\ef,
\ee
(this is possible because repeated indices give zero, from the
nilpotence of $\psi_i$ variables).
Reordering the $\psi_i$'s next to the factors with the corresponding
indices, and performing the sum over indices $i_{\alpha}$'s, gives
\be
\label{eq.6535654btausym5}
\frac{1}{n!}
\int 
\mathcal{D}(\psi, \psibar)
\big( \psibar ( U + A (n-1) ) \psi \big)
\big( \psibar ( U + A (n-2) ) \psi \big)
\cdots
\big( \psibar U \psi \big)
\ef,
\ee
which coincides with (\ref{eq.43254464}), as was to be proven.
\qed

\noindent
Our case of interest is recovered by setting $U=XY$ and $V=A$. The
hypothesis (\ref{eq.789765454}) holds, as a consequence of Lemma
\ref{commute} specialized to $B=I$ (of which, because of Lemma
\ref{lem.XABcond}, the hypotheses are satisfied), while the hypothesis
(\ref{eq.symcdet91}) is verified by observing that, for any
permutation $\tau$, the three matrices $X$, $Y^{\tau}$ and $A^{\tau}$
satisfy the hypotheses of Proposition \ref{theo.capCB}(a), and by
applying (\ref{id.coldet}) to the \emph{left} hand side of the proposition
statement (we use at this aim the fact that $Y$ has weakly
row-symmetric commutators, as implied by the hypotheses of Proposition
\ref{prop:grass}(a)). Conversely, equations
(\ref{id.coldet}) and (\ref{eq.symcdet91}) are \emph{not} immediately
related, as, because of the factors $n-j$ in $Q^{\rm col}$, the matrix
on the left hand side of (\ref{eq.symcdet91}) does not correspond to
the action of $\tau$ from the right.

Remark that, with respect to Proposition \ref{prop.old}, the
level of generality of this proposition in comparison to the
specialization pertinent to Capelli-like identities is less
pronounced. This is mainly due to the fact that the hypothesis
(\ref{eq.symcdet91}) is in fact very demanding. Indeed, it implies in
particular that, for any permutation $\tau$ and any transposition
$(j\,j+1)$ of consecutive elements,
\be
\coldet (U^{\tau} + Q^{\rm col}(V^{\tau})) 
+
\coldet (U^{\tau \circ (j\,j+1)} + Q^{\rm col}(V^{\tau \circ (j\,j+1)})) 
=0
\ef.
\ee
Using the representation (\ref{coldetint}) of column-determinants,
gives
\be
\label{eq.41434354bvvcnbvc}
\int d \psi_n \cdots d \psi_1
\;
L
\,
\big[
\big( \psi ( U + V (n-j) ) \big)_r
\big( \psi ( U + V (n-j-1) ) \big)_s
+
(r \leftrightarrow s)
\big]
\,
R
=0
\ef,
\ee
where $L$ and $R$ are appropriate factors, corresponding to the
product of $(\psi(U+Q^{\rm col}))_i$ for $i \neq j, j+1$.  A
sufficient condition for the integral to vanish is that the
combination in square brackets is zero. Strictly speaking, this is not
also necessary, but it is hard to imagine a different mechanism for
the quantity above to vanish, and still the original
column-determinant being non-trivial. So we keep on investigating
under which conditions on $U$ and $V$ we have, for every $r$, $s$ and $j$,
\be
\big( \psi ( U + V (n-j) ) \big)_r
\big( \psi ( U + V (n-j-1) ) \big)_s
+
(r \leftrightarrow s)
=0
\ef.
\ee
Matching the terms with different degree in $j$ gives
\begin{align}
\{ (\psi V)_r, (\psi V)_s \} 
&= 0
\ef;
\\
\label{eq.543654344}
\{ (\psi U)_r, (\psi V)_s \} 
&=
-
\{ (\psi U)_s, (\psi V)_r \} 
\ef;
\\
\label{eq.543654344c}
\{ (\psi U)_r, (\psi U)_s \} 
&=
(\psi U)_r (\psi V)_s
+
(\psi U)_s (\psi V)_r
\ef.
\end{align}
Incidentally, equation (\ref{eq.543654344}) implies
(\ref{eq.789765454}), thus the three equations above are sufficient
for Proposition \ref{prop.comeP13perGrass} to apply.

\section*{Acknowledgments}

This work is a continuation of a previous work, done in
collaboration with A.\,D.~Sokal. As always in these cases, it is hard to
``trace a boundary'' on authorship. It is clear that the many
discussions together had a prominent role in the genesis of the
present paper.

We thank the Isaac Newton Institute for Mathematical Sciences,
University of Cambridge, for support during the programme on
Combinatorics and Statistical Mechanics (January-June 2008), where a
large fraction of this work has been done.

S.C. thanks the Universit\`{e} Paris Nord for the support offered to visit LIPN where this work 
has been finished.



\begin{thebibliography}{99}

\bibitem{uscapelli}
S.~Caracciolo, A.\,D.~Sokal and A.~Sportiello,
{\it Noncommutative determinants, Cauchy--Binet formulae, and
  Capelli-type identities.
I.~Generalizations of the Capelli and Turnbull identities,}
{\em Electron. J. Combin.}\/ {\bf 16(1)}, \#R103 (2009),\\
{\tt arXiv:0809.3516} 


\bibitem{Capelli_1882}  A.~Capelli, 
{\it Fondamenti di una teoria generale delle forme algebriche,}
Atti Reale Accad.\ Lincei, Mem.\ Classe Sci.\ Fis.\ Mat.\ Nat.\ (serie
3) {\bf 12}, 529--598 (1882)

\bibitem{Capelli_1887}  
A.~Capelli, 
{\it Ueber die Zur\"uckf\"uhrung der Cayley'schen Operation $\Omega$
  auf gew\"ohnliche Polar-Operationen,}
Math.\ Annalen\ {\bf 29}, 331--338 (1887)

\bibitem{Capelli_1888}  
A.~Capelli, 
{\it Ricerca delle operazioni invariantive fra pi\`u serie di
  variabili permutabili con ogni altra operazione invariantiva fra le
  stesse serie,}
Atti Reale Accad.\ Sci.\ Fis.\ Mat.\ Napoli (serie 2) {\bf 1}, 1--17
(1888)

\bibitem{Capelli_1890}  
A.~Capelli, 
{\it Sur les op\'erations dans la th\'eorie des formes alg\'ebriques,}
Math.\ Annalen {\bf 37}, 1--37 (1890)

\bibitem{Weyl_46}
H.~Weyl, 
{\it The Classical Groups, Their Invariants and Representations,} 2nd
ed., Princeton Univ.\ Press,
1946

\bibitem{Procesi}
C.~Procesi, 
{\em Lie Groups: An Approach through Invariants and Representations,}
Springer Verlag, 2007

\bibitem{Howe_89}
R.\ Howe,
{\it Remarks on classical invariant theory,}
Trans.\ Amer.\ Math.\ Soc.\ {\bf 313}, 539--570 (1989),
and erratum {\bf 318}, 823 (1990).


\bibitem{Turnbull_48}
H.W.\ Turnbull,
{\it Symmetric determinants and the Cayley and Capelli operators,}
Proc.\ Edinburgh Math.\ Soc.\ {\bf 8}, 76--86 (1948).

\bibitem{Wallace_53}
A.H.\ Wallace,
{\it A note on the Capelli operators associated with a symmetric
  matrix,}
Proc.\ Edinburgh Math.\ Soc.\ {\bf 9}, 7--12 (1953).

\bibitem{Kostant_91}
B.\ Kostant and S.\ Sahi,
{\it The Capelli identity, tube domains, and the generalized Laplace
  transform,}
Adv.\ Math.\ {\bf 87}, 71--92 (1991).

\bibitem{Kostant_93}
B.\ Kostant and S.\ Sahi,
{\it Jordan algebras and Capelli identities,}
Invent.\ Math.\ {\bf 112}, 657--664 (1993).


\bibitem{Howe_91}
R.\ Howe and T.\ Umeda,
{\it The Capelli identity, the double commutant theorem, and
  multiplicity-free actions,}
Math.\ Ann.\ {\bf 290}, 565--619 (1991).

\bibitem{Itoh_00}
M.\ Itoh,
{\it Capelli elements for the orthogonal Lie algebras,}
J.\ Lie Theory {\bf 10}, 463--489 (2000).

\bibitem{Itoh-Umeda_01}
M.\ Itoh and T.\ Umeda,
{\it On central elements in the universal enveloping algebras of the
  orthogonal Lie algebras,}
Compositio Math.\ {\bf 127}, 333--359 (2001).

\bibitem{Noumi_94}
M.\ Noumi, T.\ Umeda and M.\ Wakayama,
{\it A quantum analogue of the Capelli identity and an elementary
  differential calculus on ${\rm GL}\sb q(n)$,}
Duke Math.\ J.\ {\bf 76}, 567--594 (1994).

\bibitem{Noumi_96}
M.\ Noumi, T.\ Umeda and M.\ Wakayama,
{\it Dual pairs, spherical harmonics and a Capelli identity in quantum
  group theory,}
Compositio Math.\ {\bf 104}, 227--277 (1996).

\bibitem{Umeda_98}
T.\ Umeda,
{\it The Capelli identities, a century after,}
Amer.\ Math.\ Soc.\ Transl.\ Ser.\ 2 {\bf 183}, 51--78 (1998).


\bibitem{Kinoshita_02}
K.\ Kinoshita and M.\ Wakayama,
{\it Explicit Capelli identities for skew symmetric matrices,}
Proc.\ Edinburgh Math.\ Soc. {\bf 45}, 449--465 (2002).

\bibitem{Mukhin_06}
E.\ Mukhin, V.\ Tarasov and A.\ Varchenko,
{\it A generalization of the Capelli identity,}
in
{\it Algebra, Arithmetic and Geometry -- Manin Festschrift, vol.\;II,}
Yu.\ Tschinkel and Yu.\ Zarhin eds., 
Progr.\ in Math.\ {\bf 270}, 383--398 (2007).
{\tt arXiv:math/0610799}

\bibitem{Molev_99}
A.\ Molev and M.\ Nazarov,
{\it Capelli identities for classical Lie algebras,}
Math.\ Ann.\ {\bf 313}, 315--357 (1999).

\bibitem{Nazarov_97}
M.\ Nazarov,
{\it Capelli identities for Lie superalgebras,}
Ann.\ Sci.\ \'Ecole Norm.\ Sup.\ {\bf 30}, 847--872 (1997).

\bibitem{Nazarov_98}
M.\ Nazarov,
{\it Yangians and Capelli identities,}
Amer.\ Math.\ Soc.\ Transl.\ Ser.\ 2 {\bf 181}, 139--163 (1998).


\bibitem{Okounkov_96b}
A.\ Okounkov,
{\it Young basis, Wick formula, and higher Capelli identities,}
Internat.\ Math.\ Res.\ Notices {\bf 17}, 817--839 (1996).










\bibitem{Manin_87}  
Yu.\,I.~Manin, 
{\it Some remarks on Koszul algebras and quantum groups,}
Ann.\ Inst.\ Fourier (Grenoble) {\bf 37},
191--205 (1987)

\bibitem{Manin_88}  
Yu.\,I.~Manin, 
{\it Quantum Groups and Non-Commutative Geometry,}
Centre de Recherches Math\'ematiques, Universit\'e de Montr\'eal, 1988
 
\bibitem{Manin_91}  
Yu.\,I.~Manin, 
{\em Topics in Noncommutative Geometry,}
Princeton Univ.\ Press,
1991

\bibitem{Chervov_08}  
A.~Chervov and G.~Falqui,
{\it Manin matrices and Talalaev's formula,}
J.\ Phys.\ A: Math.\ Theor.\ {\bf 41}, 194006 (2008)
{\tt arXiv:0711.2236}

\bibitem{Chervov_09}
A.~Chervov, G.~Falqui and V.~Rubtsov,
{\it Algebraic properties of Manin matrices 1,}
Adv.\ Appl.\ Math.\ {\bf 43} 239--315 (2009) 
{\tt arXiv:0901.0235}

\bibitem{KonvPhd}
M.~Konvalinka,
{\it Combinatorics of determinantal identities,}
PhD Thesis, MIT, June 2008 (Supervisor: I.~Pak),
{\tt http://math.mit.edu/$\sim$konvalinka/thesis.pdf} or
{\tt http://www.fmf.uni-lj.si/$\sim$konvalinka/thesis.pdf}

\bibitem{Flajolet}
P.~Blasiak and P.~Flajolet,
{\it Combinatorial Models of Creation-Annihilation,}\\
{\tt arXiv:1010.0354}

%
%
%
%
%
%

\bibitem{saglewalde}
A.\,A.~Sagle and R.\,E.~Walde,
{\it Introduction to Lie Groups and Lie Algebras}, 
Academic Press,
1973

\bibitem{Dirac}
P.\,A.\,M.~Dirac,
{\it The Principles of Quantum Mechanics}, 
4th ed., Oxford Univ.\ Press, 1958

\bibitem{uscayley}
S.~Caracciolo, A.\,D.~Sokal and A.~Sportiello,
{\em Algebraic/combinatorial proofs of Cayley-type identities for derivatives of
  determinants and pfaffians}, 
  Adv.\ Appl.\ Math.\ {\bf 50} 479--594 (2013) 
{\tt arXiv:1105.6270}


\bibitem{Flajolet-book}
P.~Flajolet and R.~Sedgewick,
{\em Analytic Combinatorics},
Cambridge Univ.\ Press, 2009

\bibitem{Cyril}
C.~Banderier and P.~Flajolet,
{\it Basic Analytic Combinatorics of Directed Lattice Paths,}
Theor.\ Comp.\ Science {\bf 28} 37-80 (2002) 

\bibitem{cbh_camp}
J.~Campbell, 
{\it On a Law of Combination of Operators bearing on the Theory of
  Continuous Transformation Groups,}
Proc.\ London Math.\ Soc.\ s1-{\bf 28} 381-390 (1896);
{\it On a Law of Combination of Operators (Second Paper),}
Proc.\ London Math.\ Soc.\ s1-{\bf 29} 14-32 (1897) 

\bibitem{cbh_poinc}
H.~Poincar\'e,
{\it Sur les Groupes Continus,}
Camb.\ Philos.\ Trans.\ {\bf 18} 220-255 (1899) \\
{\small {\tt http://www.archive.org/download/transactions18camb/transactions18camb.pdf}}

\bibitem{cbh_baker}
H.~Baker, 
{\it Further Applications of Metrix Notation to Integration Problems,}
Proc.\ London Math.\ Soc.\ s1-{\bf 34} 347-360 (1901);
{\it On the Integration of Linear Differential Equations,}
Proc.\ London Math.\ Soc.\ s1-{\bf 35} 333-378 (1902);
{\it Alternants and Continuous Groups,}
Proc.\ London Math.\ Soc.\ s2-{\bf 3} 24-47 (1905) 

\bibitem{cbh_haus}
F.~Hausdorff,
{\it Die symbolische Exponentialformel in der Gruppentheorie,}
Ber.\ Verh.\ Saechs.\ Akad.\ Wiss.\ Leipzig {\bf 58} 19-48 (1906) 

\end{thebibliography}
\end{document}